\documentclass[11pt]{article}
 \usepackage{amsmath, amsthm, amsfonts, amssymb}

\hoffset= - 1.1 cm \voffset=  0.0 cm \numberwithin{equation}{section}

\newtheorem{theorem}{Theorem}[section]
\newtheorem{lemma}[theorem]{Lemma}
\newtheorem{proposition}[theorem]{Proposition}
\newtheorem{corollary}[theorem]{Corollary}
\newtheorem{remark}[theorem]{Remark}
\newtheorem{definition}[theorem]{Definition}
\newtheorem{example}[theorem]{Examples}
\newtheorem{hypothesis}[theorem]{Hypothesis}

\def\qed{\hfill \hbox{\hskip 6pt\vrule
width6pt height6pt depth1pt  \hskip1pt}
\smallskip}

\newcommand{\bth}{\begin{theorem}}
\renewcommand{\eth}{\end{theorem}}
\newcommand{\bpr}{\begin{proposition}}
\newcommand{\epr}{\end{proposition}}
\newcommand{\bco}{\begin{corollary}}
\newcommand{\eco}{\end{corollary}}
\newcommand{\ble}{\begin{lemma}}
\newcommand{\ele}{\end{lemma}}
\newcommand{\bpf}{\begin{proof}}
\newcommand{\epf}{\end{proof}}
\newcommand{\bex}{\begin{example}}
\newcommand{\eex}{\end{example}}
\newcommand{\bdf}{\begin{definition}}
\newcommand{\edf}{\end{definition}}
\newcommand{\bre}{\begin{remark}}
\newcommand{\ere}{\end{remark}}
\def\ds{\displaystyle}

\newcommand{\beq}{\begin{equation}}
\newcommand{\eeq}{\end{equation}}
\newcommand{\bal} {\begin{aligned}}
\newcommand{\eal}{\end{aligned}}
\newcommand{\ben}{\begin{enumerate}}
\newcommand{\een}{\end{enumerate}}
\newcommand{\beqr}{\begin{eqnarray*}}
\newcommand{\eeqr}{\end{eqnarray*}}
\def\aa{\tilde{\cal A}}

\newcommand{\bit}{\begin{itemize}}
\newcommand{\eit}{\end{itemize}}

\def\P{{\mathbb P}}
\def\R{{\mathbb R}}

\def\E{{\mathbb E }}
\def\Q{{\mathbb Q}}
\def\Z{{\mathbb Z}}
\def\N{{\mathbb N}}

\def\lan{{\langle}}
\def\ran{{\rangle}}
\def\hh{{\vskip 1mm} \noindent}
\def\vv{{\vskip 1mm}}

\def\nor{\,    {  | \, \!   | \,  } }

\begin{document}


\title {\Huge { Global
Schauder estimates for  a class of degenerate Kolmogorov equations
 }}

\maketitle

\begin{center}

 Enrico Priola \footnote
 {Partially supported by  the Italian National Project
 MURST ``Equazioni di Kolmogorov''.}

\vspace{ 2 mm } {\small  \it Dipartimento di Matematica,
Universit\`a di Torino, \par
  via Carlo Alberto 10,  \ 10123, \ Torino, Italy. \par
 e-mail \ \ enrico.priola@unito.it }
\par \ \par

\end{center}

{\vskip 7mm }
 \noindent {\bf Mathematics  Subject Classification (2000):} \
 \  35K65, 35J70, 47D07, 35B65.

 \ \par

\noindent {\bf Key words:} Schauder estimates, degenerate elliptic
and parabolic equations,  diffusion Markov semigroups.

\vspace{2.5 mm}

\noindent {\bf Abstract:} We consider a class of possibly
degenerate second order elliptic operators $\cal
 A$ on $\R^n$. This class includes  hypoelliptic
 Ornstein-Uhlenbeck type operators having  an
  additional  first order term with unbounded coefficients.
    We establish global  Schauder estimates in H\"older spaces
     both for  elliptic
  equations  and for parabolic Cauchy problems involving ${\cal A}$.
  The H\"older function spaces are defined with respect to a
  non-euclidean metric related to the operator $\cal A$.

\section {Introduction}

 Let us consider the following  possibly degenerate
 second order elliptic operator $\cal A$ on $\R^n$:
 \beq \label{op} \bal
 {\cal A} u (x) &=
  \mbox{$\frac{1}{2}$} {\mbox{\rm Tr }}  (Q D^2 u (x)) +
  {\langle}  Ax, D u (x) {\rangle}
 + {\langle}  F(x), D u (x){\rangle}
 \\ &= {\cal A}_0 u (x) + {\langle}  F(x), D u (x){\rangle},
 \;\;\;\; x \in \R^n.
\eal \eeq
  Here  $Q$ and  $A$ are $n \times n$ real matrices,
$Q$ is
  symmetric  and non-negative definite, Tr$(\cdot) $ denotes the
  trace and ${\langle}  \cdot, \cdot {\rangle} $ the   inner product in
  $\R^n$. Moreover
    $F : \R^n \to \R^n$ is a possibly unbounded regular vector field.
  Degenerate Kolmogorov operators like $\cal A$ arise in
  Kinetic Theory and in  Mathematical Finance (see, for instance,
 \cite{DV}, \cite{DP} and the references therein).
    Moreover, the operator $\cal A$ contains in the special case
     of $F=0 $ the well-studied possibly degenerate
     Ornstein-Uhlenbeck operator ${\cal A}_0$.

      The aim of this paper is to prove global  Schauder estimates for
 elliptic equations and
 parabolic Cauchy problems involving
  the operator $\cal A$.   We obtain
   optimal regularity results in  H\"older spaces  for both
 \beq  \label{1}
  \lambda u (x) - {\cal A} u(x) =  f (x), \;\;\; x \in \R^n,\; \;\;
  \mbox{and}
\eeq \beq \label{2}  \left\{ \bal \partial_t v(t,x) & = {\cal A} v
(t,x) \, +\, H(t,x), \;\;\; t \in (0,T],
 \; x \in \R^n,\\
      v(0, x) & = g(x),\;\; x \in \R^n,
  \eal \right.
\eeq
 where  $\lambda >0$ and the  functions  $f$, $g$ and $H$ are
 assigned.  Let us collect our assumptions on the operator $\cal A$
 (compare with \cite{Pr}).
 \begin{hypothesis} \label{hy} {\em { (i)} \
  the  symmetric matrix
  $Q = (q_{ij})_{i,j=1, \ldots ,n}$
   is given by
   \beq \label{q0}
Q =  \left ( \begin{matrix} Q_0 & 0 \\ 0 & 0
  \end{matrix} \right ), \; \mbox{where $Q_0$ is a
   positive definite $\tilde p \times \tilde
   p\,$-matrix}, \, \, 1 \le {\tilde p} \le n;
  \eeq
  $\nu_1$ and $\nu_2$ stand  for
  the  smallest and the largest
   eigenvalue of $Q_0$ respectively ($0 < \nu_1 \le \nu_2$);

\hh  { (ii) } \ the vector field $F: \R^n \to \R^n$  has the form
$F(x)= (F_1(x), \ldots, F_{\tilde p}(x), 0, \ldots, 0)$, $x \in
 \R^n$, i.e., $F(x) \in $Im(Q), for any $x \in \R^n$;

\hh { (iii)} the non-zero  coefficients of $F$, $F_i : \R^n \to
\R$, $i=1, \ldots , {\tilde p},$ are Lipschitz
 continuous functions having bounded partial
 derivatives up to the third
 order on $\R^n$;

\hh { (iv)} \ there exists a  nonnegative integer $k$,
 such
 that the vectors
\beq \label{kal}
 \{ e_1 , \ldots, e_{\tilde p}, A e_1 ,
 \ldots, Ae_{\tilde p}, \ldots, A^k e_1 , \ldots, A^k
 e_{\tilde p}\} \;\;\; \mbox{generate} \;\; \R^n
 \eeq
($e_1 , \ldots, e_{\tilde p}$ are the first
 ${\tilde p}$ elements of the canonical basis in $\R^n$); we
 denote by $k$ the {\it smallest} nonnegative integer such that \eqref{kal}
 holds (one has $0 \le k \le n-1$). }
\end{hypothesis}
\noindent Condition \eqref{kal} can be also written as
 Rank$[Q^{1/2}, A Q^{1/2},$ $\ldots ,$ $ A^{k} Q^{1/2} ]=n$.
 By the  well-known
   H\"ormander condition on commutators, \eqref{kal}
      is equivalent to the   hypoellipticity of the operator
   ${\cal A}_0 - \partial_t $
  in $(n+1) $  variables  $(t, x_1, \ldots,  x_n)$;
   see \cite{LP}.
   Our  operator $\cal A$
 has the following expression
 $$
   {\cal A} u (x) = \frac{1}{2} \sum_{i,j =1}^{{\tilde p}} q_{ij}
\partial_{x_i x_j}^2
   u (x) +
  \sum_{i =1}^{{\tilde p}} F_{i}(x) \partial_{x_i } u (x) +
  \sum_{i,j =1}^n a_{ij} \, x_{j} \partial_{x_i } u (x), \;\;\;
  x \in \R^n,
 $$
 where  the  $a_{ij}$ are the components of the
 matrix $A$ and $\partial_{x_i}$ and
 $ \partial_{x_i x_j}^2$ are  partial derivatives.
  Clearly, the operator $\cal A$ is non-degenerate only when $\tilde p =
  n$ (this implies $k=0$).

  Let us explain Schauder estimates for
  \eqref{1} and \eqref{2}.
 In  the elliptic equation \eqref{1} we assume that  $f \in {\cal C}^{\theta}_d
 (\R^n)$, $\theta \in (0,1)$, i.e., $f$ is a  real bounded function
 on $\R^n$, which is H\"older continuous with respect to a suitable
 non-euclidean metric $d$ related to $\cal A$.
 We show  that \eqref{1} has a unique bounded distributional
  solution $u
\in {\cal C}^{2 + \theta}_d (\R^n)$,  and that
 there exists $C
>0$, independent of $f$ and  $u$, such that $\| u\|_{2 + \theta, d} \le
 C \, \| f\|_{\theta, d}$. Note that this implies
 $$
 \mbox{$ \| u\|_{0} +  \sum_{i,j=1}^{\tilde p}\| \partial^2 _{x_i x_j}
  u \|_{\theta , d}  \, \le \, C \, \| f\|_{\theta, d},
 $}$$
 where $\| u\|_0$ denotes the sup-norm of $u$
 (see Theorem \ref{buon}).
   Concerning the  Cauchy problem \eqref{2} we prove analogous
    parabolic  Schauder  estimates,
  assuming that  $g \in  {\cal C}^{2 + \theta}_d (\R^n)$ and
  $H (t, \cdot ) \in {\cal C}^{\theta}_d (\R^n)$, uniformly in
  $t \in [0,T]$ (see Theorem \ref{buon1}).
   We refer to   Section 2 for a precise definition of the metric $d$.
   Here we  give an   example of  $d$.  We  consider
   the following two-dimensional operator $\cal
  A$,
 \beq \label{exa}
 {\cal A} u (x,y) = \mbox{$\frac{1}{2}$} \partial_{xx}^2 u (x,y) \, +
 \,
   F_1(x,y) \partial_{x}
 u (x,y)
 \,  +  \, (x + y)  \partial_y u (x,y), \;\; \;\;\; (x,y)  \in
 \R^2
 \eeq
 (this operator verifies Hypothesis \ref{hy} with $\tilde p=1$ and
  $k=1$).
 In this case, the metric $d$ is given by  $d (z,z')= |x - x'|+ |y
- y'|^{1/3}$, for any $z = (x,y)$ and $z' = (x', y')
  \in \R^2$. Remark that
    $d$ is mentioned in   \cite[page 11]{St}
  and
  it is related to  certain   distances
    associated to degenerate operators
  such as
  Hormander's sum of squares of vector
  fields (see  in particular the metric
   $\rho_3$ in \cite[page 112]{NSW}).
  Moreover, $d$     is
     a special case of the parabolic  pseudo-metric
 considered in \cite{DP} (see also
  \cite{LP}).


Let's now examine   related papers on Schauder estimates. A
general
 theory of  {\it local}  regularity in Sobolev and H\"older spaces
  is available for  degenerate operators which are sum of
 squares of vector fields  (see in particular \cite{Fo}, \cite{RS} and
 \cite{Kr}).
   Local  $C^{\theta}$-estimates for operators more general than $\cal
   A$,
   in which also $q_{ij}$ are variables and time-dependent, can be
   found  in   \cite{M},
    \cite{P} and \cite {DP} (see also  the references therein).
  Concerning {\it global } regularity results  for solutions of
possibly degenerate equations like \eqref{1} and \eqref{2} in
spaces of continuous functions, we mention \cite{Lu}, \cite{Lo},
\cite{Lo1}, \cite{P}, \cite{Sa}. In \cite{Lu} Schauder estimates
are established for the  Ornstein-Uhlebeck
   operator ${\cal A}_0$
 only assuming \eqref{kal}.
 In  \cite{Lo} and \cite{Lo1} Schauder estimates are proved
  for
 Ornstein-Uhlenbeck types operators ${\cal A}_0$ when
  $F_i =0$ but
 $q_{ij}$ are not constant and  can be
  unbounded;  in  \cite{Lo} and \cite{Lo1}
   it is assumed   $k \le 1$ in hypothesis \eqref{kal}.
    Uniform estimates
   for solutions to the Cauchy problem \eqref{2}
     involving $\cal A$
    with $H=0$ are given in   \cite{Pr};
   these  are proved without any restriction on
    $k$ and are
   preliminary to
     the
     Schauder estimates
    of the present paper.
  In \cite{Sa}  Schauder
 estimates  are proved for $\cal A$ assuming $k\le 1$ in \eqref{kal}
  and  imposing   an additional hypothesis  (which is not satisfied
  in \eqref{exa}).

 To prove elliptic Schauder estimates,  one  considers
the function \beq
 \label{lap1}
 u (x) = \int_0^{+ \infty} e^{- \lambda  t} \, P_t f(x) dt,\;\; x \in
 \R^n,
\eeq where $P_t$ is the diffusion Markov semigroup associated to
${\cal A}$ (i.e., $v(t,x)= (P_t f)(x) =
P_t f(x)$ provides the
classical
 solution to \eqref{2} when $H=0$, see   \cite{Pr}). The function
  $u$ is the unique bounded distributional solution to
 \eqref{1} (see Theorem \ref{eli}).
 One
 proves global regularity properties for $u$ by means of
  sharp $L^{\infty}$-estimates on the
 spatial partial derivatives of $P_t f$ involving the H\"older norm
  of $f$ (the behaviour in $t$ of such estimates
    as $t$ tends to
  $0^+$ is crucial).
  This is the basic idea
 indicated in \cite{DL}
 in order to study Schauder estimates for non-degenerate Kolmogorov
   operators.
  This method  has been  much used in  recent papers
    also in
  combination  with  \cite{Lu1}
  (see  \cite[Chapter 1]{Ce}, \cite[Chapter 6]{BL} and
   the references therein).
  In \cite{Lu} the  $L^{\infty}$-estimates have
  been proved  using the  explicit formula   of
  the Ornstein-Uhlenbeck semigroup $P_t$ associated to ${\cal A}_0$
  (which is  not available
  when $F \not =0$ in ${\cal A}$). In \cite{Lo}, \cite{Lo1} and \cite{Sa}
   the uniform  estimates
   are obtained by
  a priori estimates of Bernstein type
   combined  with
    an interpolation result proved in
     \cite[Lemma 5.1]{Lo1} when $k\le 1$.
  We get  the   $L^{\infty}$-estimates
 involving  H\"older norms
    in Theorem \ref{sti} by working directly
  on  some probabilistic
 formulae for  the spatial derivatives of $P_t f $
  (which replace the explicit formulae used in \cite{Lu}).
 Such  formulae
  have been obtained in
  \cite{Pr} using Malliavin Calculus (see also \cite{B}, \cite{KS}
   and \cite{Fu}).

We believe that the probabilistic
    approach used to derive $L^{\infty}$-estimates
 could be useful in other situations.
 In particular, we  have in mind
  degenerate Kolmogorov  operators $\cal A $  in which the
 drift vector field $Ax +  F$ is replaced by
  a $C^{\infty}\,$-vector field $G : \R^n \to \R^n$;
  one  assumes that  $G$ has  all bounded derivatives and that
     there exists an integer $k$ such that
  $ e_1, \ldots, e_{\tilde p} $ and
   $G$  together with  their commutators of length at most
 $k$ span $\R^n$ at each point $x \in \R^n$.
 This problem is largely open.

  Once the previous $L^{\infty} $-estimates  are proved
   for a class of Kolmogorov operators,
   recent papers use an interpolation result
  of \cite{Lu1}
  in order to obtain Schauder estimates  for $u$ (see, for instance, \cite{Lu},
  \cite[Chapter 1]{Ce},  \cite{Lo1}, \cite{BL}, \cite{Sa}).
  We propose in Theorem \ref{buon} a
    direct approach
 to get elliptic Schauder estimates
   (this method
   applies also to parabolic Schauder estimates).

  In order to study the parabolic Cauchy problem one proceeds initially
 as in the
 elliptic case, replacing the formula \eqref{lap1} with the
 variation of constant  formula (see \eqref{cau}).
  However,  the  parabolic Schauder estimates
   are more difficult to prove than the corresponding elliptic
   ones (see  Remark
 \ref{elli}).
 In particular, they require the
  hard  estimate
  $
 \| (P_t  g) (\cdot)\|_{2 + \theta, d} $ $ \le  C \| g\|_{2 + \theta, d},
 $
 for any  $g \in C^{2 + \theta}_d (\R^n)$, $t \ge 0$, where $C $ is
 independent of $t \, \mbox{and} \, g$.

\vskip 1mm After some preliminaries contained in Section 2,  in
Section 3 we  prove
  the  $L^{\infty}$-estimates for the spatial derivatives of
  $P_t f$ involving the H\"older norm of $f$.
  In Section 4 we show that \eqref{1}  has a  unique
  distributional solution and prove elliptic Schauder estimates using
  the results of Section 3. We also establish  existence and uniqueness
  of  space-distributional solutions to the parabolic Cauchy problem
  \eqref{2} and prove the parabolic Schauder estimates.
  In the final part of the paper we
 consider more
general operators $\tilde {\cal A}$ with variable coefficients
 $q_{ij}(x)$.
   We require that the matrix $Q(x)$ has the form \eqref{q0} where
  the ${\tilde p} \times {\tilde p}$
    matrix $Q_0 (x)$ is uniformly positive; moreover, we assume  that
 $q_{ij} $ are $\theta$-H\"older continuous and that there exists
   $  \lim_{x \to \infty} Q_0(x) = Q^{\infty}_0$.
  We  obtain   elliptic and parabolic Schauder estimates
 for $\tilde {\cal A}$,   using a well known
  method based on maximum principle, a priori
estimates and continuity method (compare with \cite{Lu}). Further
extensions  of our results are proposed in Remark \ref{fine}.

{\vskip 1mm} We will use the letter $c$ or $C$ with subscripts for
finite positive
 constants whose  precise  value is unimportant;  the constants   may
 change from
    proposition to proposition.



\section{Preliminaries and notation}

 We denote by $|\cdot |$ and $\lan \cdot , \cdot \ran$ the
 euclidean norm and the standard inner product in $\R^n$ and by $\|
\cdot \|_{L}$ the operator norm in  the Banach space $L (\R^n)$ of
real $n \times n$ matrices.  If $X$ and $Y$ are real Banach
spaces, $L (X, Y)$
 denotes the Banach space of all bounded and linear operators from
  $X$ into $Y$ endowed with the operator norm.

  Let $G :  \R^n \to \R^m$ be a   mapping.
   We  denote by $DG (x)$, $D^2G(x)$ and $D^3G(x)$ respectively
 the first, second and third Fr\'echet derivative of $G$ at $x \in \R^n$
 when they exist (if  $G$ also depends on $t$,
   we write $D_x G(t,x)$, $D^2_{xx} G(t,x)$ and $D^3_{xxx} G(t,x)$).
 We have  $DG(x)[u]$, $D^2G(x) [u][v]$
  and  $D^3G(x) [u][v][w] \in \R^m$, for $u,v,w \in \R^n$.
 If $G$ is bounded, we set  $\| G \|_0 = \sup_{x \in \R^n}
 |G(x)|_{\R^m}$.

 Recall that   hypothesis \eqref{kal} is known as the Kalman
 condition   in control theory (see \cite{Za}).
  It
   is also equivalent to requiring that the following symmetric matrix
$Q_t$,
 \beq \label{qtt}
  Q_t = \int_0^t e^{s A^*} Q e^{sA} ds
 \eeq
 is positive definite for any $t>0$ (here $e^{sA}$ denotes the exponential
  matrix of $A$
 and $A^*$ the adjoint matrix of $A$).

 As in \cite{Lu} we define  an {\it orthogonal decomposition}
of $\R^n$ related to the Kalman condition \eqref{kal}.  We
consider the first ${\tilde p}$ elements $\{ e_1, \ldots,
e_{\tilde p} \}$ of the
 canonical basis in $\R^n$, $1 \le {\tilde p}
  \le n$, and introduce the subspace $V_0 =$  Span$\{e_1, \ldots
 , e_{\tilde p} \}$. Then  set
  $V_m = Im Q^{1/2} + $ $... + Im \, (A^m Q^{1/2})$
 $=$ Span$\{ e_1, \ldots
 , e_{\tilde p}, Ae_1,$ $ \ldots, A e_{\tilde p}, \ldots $ $A^m e_1,
  \ldots , A^{m}e_{\tilde p}
  \}$, for $1 \le m \le k $.
 One has $V_m \subset V_{m+1}$ and $V_k = \R^n$.  Let $ W_0 = V_0,
$  $W_1$ be the orthogonal complement of $V_0 $ in $V_1$, $W_m $
be the orthogonal complement of $V_{m-1} $ in $V_m $, for $1 \le m
\le k$. Defining the orthogonal projections $E_m$ from $\R^n $
onto $W_m$, one
 has $E_m (\R^n) = W_m$ and
  \beq\label{sp} \mbox{ $
\R^n = \bigoplus_{m=0}^k E_m (\R^n),
 $ } \eeq
 We  complete  $\{ e_1, \ldots, e_{\tilde p} \}$ in order to
  get a {\it reference orthonormal basis} $  \{ e_i \}_{i=1, \ldots ,
  n}$
   in $\R^n$ related to \eqref{sp}. This consists of generators of
    the subspaces
 $E_m (\R^n)$, $0 \le m \le k$,
  and  will be used throughout the paper.  Note that,  writing
  the operator $\cal A $ in the   coordinates
  associated to  the new basis
  the second order term Tr$(Q D^2)$ does not change.
   In the sequel $D_i$, $D_{ij}^2 $,
    and
   $D^3_{ijr}$ will
denote respectively first, second and third partial derivatives
 with respect to $\{ e_i\}$
(one
   can assume that $\{ e_i \}$ is the canonical basis  if
 $k\le 1$, compare with \cite{Lo1} and \cite{Sa}).
  Define $I_m$ as the set of indices $i$ such that $e_i$ spans
 $E_m(\R^n)$, $0 \le m \le k$. We have
 $$
 I_0 = \{ 1, \dots , {\tilde p}\}.
 $$
The metric $d$  associated to the  operator $\cal A$ is defined
using the decomposition \eqref{sp}.
 One first introduces the   {\it
quasi-norm} $\nor \cdot \nor $,
  $  \nor  x \nor := \sum_{h=0}^{k} |E_h x|^{1/ (2h+1)},$ $  x \in \R^n.
$
 Then we set
 \beq \label{di}
 d (x, y) :=  \nor x- y \nor = \sum_{h=0}^{k} |E_h
(x-y)|^{\frac{1} { 2h+1}}
 ,\;\;\; x, \, y \in \R^n.
  \eeq
 Let us introduce some  function spaces. First we consider
  euclidean function spaces and then function  spaces
 related to the metric $d$.

{\vskip 1mm} We denote by $B_b (\R^n)$ the Banach space
 of  all Borel and bounded
   functions \ $f: \R^n \to \R
$, endowed with the supremum  norm $\| \cdot\|_{0} $;
  $C_b
(\R^n)$ is the closed subspace of $B_b (\R^n)$
 consisting of all  {\it uniformly continuous and bounded}
   functions.

  $ C_b^{j}(\R^n) $, $j \in \Z_+$, $ j \ge 1$,  is
   the Banach space  of
  all $j$-times differentiable functions $f : \R^n \to \R$,
    whose partial  derivatives,
  $D_{\alpha} f$, $ \alpha \in \Z_+ ^n,$
   are uniformly continuous and bounded on $\R^n$ up to order $j$. This  is
    a  Banach space endowed with the norm $\| \cdot \|_j$,  $
 \| f {\|}_{j}\;=\;\| f{\|}_0 \;+\;\sum _{ |\alpha | \le j}^{} \;\;\|
 D_{\alpha} f{\|}_0,\;\; f \in C_b^j(\R^n).
  $ We set $C_{b}^{\infty} (\R^n) = \cap_{j \ge 1} C_b^{j}(\R^n) $.
 Moreover $C_0^{\infty}
 (\R^n) $ is the space of all functions $f \in C^{\infty}_b
 (\R^n)$ having compact support.

{\vskip 2mm} Fix $\theta \in (0,1)$.  The space  $C^{\theta}_b
(\R^n)$   stands for the  Banach space of all $\theta-$H\"older
continuous and bounded functions on $\R^n$ endowed with the  norm
$\| \cdot \|_{\theta}$, i.e., $\| f \|_{\theta } \, = $ $ \| f
\|_{0}  \,$ $ + \, [ f ]_{\theta },$ \ $f \in  {
C}_b^{\theta}(\R^n ), $ where $ [f]_{\theta } \,= \;  \sup _{z,w
\in \R^n, \; z \not = w } $ $\; \frac { | f(z)\, - \, f(w) | } { |
z-w|^{ \theta} } \; < \; \infty.$ Moreover $C^{2+ \theta}_b (\R^n)
= \{ f \in C^2_b (\R^n)\; :\; D_{ij}^2 f \in C^{\theta}_b
(\R^n),\;\; i,j =1, \ldots ,n \}$; it is a Banach space endowed
 with the norm $\| \cdot \|_{2+ \theta} $, $
 \| f\|_{2+ \theta} = \| f\|_2  + \sum_{i,j =1}^n \| D_{ij}^2 f
 \|_{\theta},\;\; f \in C^{2+ \theta}_b
(\R^n). $  In a similar way one defines the Banach space $C^{1+
\theta}_b (\R^n)$. Next,   we define  function spaces related to
 the metric $d$.

{\vskip 1mm}
Let $\gamma \in (0,3)$ and {\it $\gamma$
 non-integer}.    We define
  ${\cal C}^{\gamma}_d (\R^n)$  as the space of all  functions $f
  \in C_b (\R^n)
 $ such that, for any $z \in \R^n $ and  for any integer $m$, $0
\le m \le k$, the map:
$$
 x \mapsto f (z + x) \;\; \mbox {\rm belongs to } \; C^{\gamma / (2 m +1)}_b \,
 (E_m
 (\R^n)),
$$
 with the  $\| f(z+ \, \cdot) \|_{\gamma/ (2m +1)}$ bounded by a
 constant {\it independent}  of $z$ (identifying
 each subspace $E_m
 (\R^n)$  with $\R^{n(m)}$, where $n(m)=$ dim$[E_m
 (\R^n)$],  the euclidean function spaces $C^{\gamma / (2 m +1)}_b \,
 (E_m
 (\R^n))$ are well defined);  ${\cal C}^{\gamma}_d (\R^n)$
  is a Banach space  with the  norm $\| \cdot\|_{\gamma,
  d}$,
 $$
\| f \|_{\gamma, d} : = \sum_{m = 0}^k  \sup_{z \in \R^n} \, \|
 f (z + \cdot) \|_{C^{\gamma / (2m +1)}_b \, (E_m (\R^n))}, \,\;\;
 f \in {\cal C}^{\gamma}_d (\R^n).
$$
It is easy to see that  if $\gamma \in (0,1)$ and $f \in
C_b(\R^n)$, then $f \in {\cal C}^{\gamma}_d (\R^n)$ if and only if
 $f$ is $\gamma-$H\"older
 continuous with respect to the metric $d$, i.e.
 \beq \label{fh}
  [f]_{\gamma, d} = \sup_{ x,\, y \in \R^n, \, x \not = y }
  \; { | f (x) - f(y)  | } {\;  \nor  x-y \nor^{- \gamma}} < +\infty.
  \eeq
 Moreover an equivalent norm in ${\cal C}^{\gamma}_d (\R^n)$, $\gamma \in
(0,1)$,  is  $ \| \cdot \|_0 + [ \, \cdot
  \, ]_{\gamma, d}.$  One can also define  $C^{\alpha}_d (\R^n)$
   for general real $\alpha >0$ (see   \cite{Lu}); we will
     only use the  spaces introduced above.



 {\vskip 1mm}
In \cite[Lemma 2.1]{Lu} it is proved that if $f \in
{\cal C}^{2 + \theta }_d (\R^n)$, $\theta \in (0,1)$, then for any
 $i, j \in I_0$, we have both $D_i f \in {\cal
 C}^{\theta +1 }_d (\R^n)$   and   $D_{ij}^2 f \in {\cal
 C}^{\theta }_d (\R^n) $; moreover there exists $C$, independent
 of $f$, such that
\beq \label{lun}
 \| D_i f \|_{1+ \theta, d} +  \|  D_{ij}^2 f \|_{\theta, d}  \le C  \,
 \| f \|_{2 + \theta, d}, \;\;\; \;\; i, j \in I_0.
 \eeq
 Let $f\in {\cal C}^{\gamma}_d (\R^n)$, $\gamma \in (2,3)$. For any $x \in \R^n$,
  we will  consider
   $D_{E_0 \,} f(x) \in \R^n$,  the gradient of
$f$ at $x \in \R^n$ in the {\it directions of $E_0 (\R^n)$,} i.e.,
 \beq \label{d0}
D_{E_0 \,} f(x) = \big( D_1 f(x), \ldots, D_{\tilde p} f (x), 0,
\ldots, 0 \big)
 \eeq
  and, similarly,
  the $n \times n$ Hessian
  matrix $  D^2_{E_0 \,} f(x) $ in the directions of $E_0 (\R^n)$,
 i.e., $\big ( D^2_{E_0 \,} f(x) \big)_{ij} = D_{ij}^2 f(x)$,
 if both  $i$ and $j \in I_0$; $\big ( D^2_{E_0 \,} f(x) \big)_{ij}=0$
 otherwise.

 \vskip 0.5 mm
    We finish the section with an equivalent
   definition of  ${\cal C}^{\gamma}_d
(\R^n)$.
   Let $f
\in C_b (\R^n)$; we introduce, for any $x,  v \in \R^n$,
 \beq \label{zyg}
\triangle_v^3 f (x) = f(x) - 3 f(x+ v) + 3 f(x+ 2v) -  f(x+ 3v).
 \eeq
\ble  \label{zyg1} Let $\gamma \in (0,3)$  non-integer.
 Let
$f \in C_b (\R^n)$. Then $f \in {\cal C}^{\gamma}_d (\R^n)$
 {\it if and only if}
$$
 [f]_{\gamma, d, 3} = \sup_{x,\, v  \, \in \R^n,\,  v \not =0,\,
   \nor v \nor \le 1}
  \;
  | \triangle_v^3 f (x)| \, { {\nor  v \nor }^{- \gamma}} \;
  \,  < \, + \infty,
$$
 see \eqref{di}.
 Moreover $\| \cdot \|_0 + [ \, \cdot \, ]_{\gamma, d, 3}$ is
equivalent to the norm $\| \cdot\|_{\gamma, d}$.
 \ele
\bpf  We use the following Triebel result (see \cite[Section
2.7.2]{T}).
 Let $g
 \in C_b (\R^n)$. Then $g$ belongs to $C^{\gamma}_b (\R^n)$, $\gamma \in
 (0,3)$ non-integer,  if and only if
   \beq \label{nag}
    [g]_{\gamma, 3} =  { \sup_{x \in \R^n, \,  |v| \le 1, \,
 v \not = 0  } } \, {  |v|^{- \gamma}} \, |\triangle_{v}^3 g
 (x)|
\, < \infty.  \eeq
 Moreover in $C^{\gamma}_b (\R^n) $ the norm $\|
\cdot \|_{\gamma}$ is equivalent to $\| \cdot \|_0 $ $+ \;  [ \,
\cdot \, ]_{\gamma, 3}$.

\hh $\Longrightarrow $ Let $f  \in {\cal C}^{\gamma}_d(\R^n)$
 and fix $v \in \R^n$. We set $v = v_0 +  v_1$, where
  $v_0 = E_0 v$ and $v_1 = \sum_{h=1}^k E_h v $ $= v - E_0 v
  $, see \eqref{sp}. We get,
   for any $x \in \R^n$,
\beqr | \triangle_{v}^3 f
 (x)|  \le  \big | f(x) - f(x + v_1)  \big |
 \\ +
\big |  f(x + v_1 ) - 3 f(x+ v_1 +  v_0) + 3 f(x+ v_1 +  2v_0) -
 f(x+ v_1 + 3v_0) \big |
\\  +
   3 \big |  f(x+ 2 v_1 +  2v_0) -  f(x+ v_1 +  2v_0) \big |
    +  \big |  f(x+ v_1 +  3v_0)  -  f(x+ 3 v_1 +  3v_0)
    \big |
\\
\le   \| f\|_{\gamma, d}  \big ( 4 \mbox{ $ \sum_{h=1}^k |E_h
v|^{\frac{ \gamma} {2h +1}} \, + \, \sum_{h=1}^k |E_h (2
v)|^{\frac{ \gamma} {2h +1}} \, + \, |v_0|^{\gamma } $} \big ) \le
C \| f\|_{\gamma, d} \nor v {\nor}^{\gamma}.
 \eeqr
 $ \Longleftarrow $ Let $f \in C_b(\R^n)$  and
  take  $v_h \in E_h (\R^n)$, with $0 \le h \le k$. By
assumption, we know that $ |\triangle_{v_h}^3 f (x)| \le
[f]_{\gamma, d,3}  |v_h|^{\gamma/ (2h +1)}$, for any $x \in \R^n$.
It follows that
 $f(x+  \cdot ) \in C^{\gamma / (2h +1)}_b \, (E_h (\R^n))$ and
 there exists $C>0$ independent of $f$ and $x$ such that
  $
\| f(x+ \, \cdot ) \|_{C^{\gamma / (2h +1)}_b \, (E_h (\R^n))}
 \le \, C \, ( \| f\|_0 +  [ f]_{\gamma, d, 3}), $ $ 0 \le h \le k.
 $ Thus $f \in {\cal C}^{\gamma}_d(\R^n).$  The proof is complete. \epf


\section{Estimates on  the diffusion semigroup associated to ${\cal A}$}

  In this section we  consider   the diffusion semigroup $P_t$
 associated to the operator $\cal A$ (compare with \eqref{lap1}).
  We   obtain $L^{\infty}$-estimates on the first, second and third
 spatial partial derivatives of $P_t f $, in terms of the H\"older-norm of
 $f$.
  These estimates
 will lead  in the next section to  Schauder estimates for \eqref{1}
   and \eqref{2}.

{\vskip 1mm} Let $(\Omega, ({\cal F}_t)_{t \ge 0}, {\cal F}, \P)$
be a complete stochastic basis (satisfying the usual assumptions;
see, for instance, \cite{IW}). Let $W_t$, $t \ge 0,$ be a standard
$n$-dimensional Wiener process defined and adapted on the
stochastic basis.
  Let $X_t^x$ be the unique (strong)  solution to the SDE
 \begin{equation} \label{xx}
X_t^x = x + \int_0^t A X_s^x ds + \int_0^t F (X_s^x) ds + Q^{1/2}
W_t, \;\;\; t \ge 0, \; x \in \R^n,
 \end{equation}
 $\P$-a.s., where the matrix $A$ is the same as in  \eqref{op} and
$Q^{1/2}$ is the unique $n \times n$ symmetric nonnegative
definite square root of $Q$.  The {\it diffusion semigroup} $P_t$
 associated to $\cal A$ is the  family of linear contractions
$P_t : B_b (\R^n) \to B_b (\R^n)$, $t \ge 0,$   defined by
 \begin{equation}  \label{pt}
 P_t g(x) := \E [g(X_t^x)],\;\;\; t \ge 0, \; g \in B_b (\R^n),\; \; x \in
 \R^n,
 \end{equation}
where the expectation is taken  with respect to $\P$. Introducing
 the Ornstein-Uhlenbeck process $Z_t^x$, which solves \eqref{xx} when $F=0$,
 \begin{equation}  \label{ou}
 Z_t^x = e^{tA} x + Z_t^0,
\;\; \text{where} \;\;
  Z_t^0 = \int_0^t e^{(t-s) A}  Q^{1/2} dW_s,
 \end{equation}
 we have:
  $ X_t^x = Z_t^x + \int_0^t e^{(t-s) A} F(X_s^x)ds.
  $

   Let us recall an application   of the Girsanov theorem
  which will be used in the proof of Theorem \ref{sti} (see also
   \cite{Pr}).  Fix  $t>0$, $x \in \R^n$,  and define
  $Q^{-1/2} =  \left ( \begin{matrix} Q_0^{-1/2} & 0 \\ 0 & 0
  \end{matrix} \right )$; then consider
 the stochastic process
 \begin{equation}  \label{ls}
L_ s^x := W_s - \int_0^s (Q^{-1/2} F) (Z_r^x) dr = W_s - \int_0^s
G (Z_r^x) dr, \;\; s \in  [0,t],
 \end{equation}
where we have set $G  := Q^{-1/2} F $. By the Girsanov theorem,
the process $L_ s^x $ is a Wiener process on $(\Omega, ({\cal
F}_s)_{s \le t}, {\cal F}_t, \Q)$, where $\Q$ is a probability
measure on $(\Omega, {\cal F}_t)$
 having  density $\Phi(t,x) $ with respect to $\P$, i.e.,
 $$
 \Q(A) := \E [1_A \, \Phi (t,x) ], \; \text{where} \;\;
  \Phi (t,x) = \exp \Big ( \int_0^t {\langle}  G(Z_s^x), dW_s {\rangle}
   \, -  \, \frac{1}{2}  \int_0^t | G(Z_s^x)
   |^2 ds  \Big),
$$ for any $A \in {\cal F}_t$.
   The processes $Z^x = (Z_s^x)$ and $X^x = (X_s^x)$, $s \in [0,t]$, satisfy the
same equation \eqref{xx}
 in $(\Omega , {\cal F}_t, \Q,
 (L_ s^x))$ and $(\Omega , {\cal F}_t, \P, (W_s)) $
 respectively. Therefore, by uniqueness,   the   laws of the
 processes $Z^x$ and $X^x$ on $C([0,t]; \R^n)$ are the same
  (under the probability
  measures
 $\Q$ and
 $\P$ respectively).
 This  implies that
 \begin{equation}  \label{33}
 \begin{aligned}
&& P_t f(x) = \E [f(X_t^x)] = \E  [f(Z_t^x )\,  \Phi (t,x)],\;\;\;
f \in B_b
 (\R^n).
\end{aligned}
 \end{equation}

\vskip 1mm  The next theorem  is proved in \cite{Pr}. It provides
 probabilistic formulae and preliminary  uniform estimates
  for the spatial partial
derivatives of $P_t f$ up to the third order (the formula for the
first derivatives was  obtained in \cite{Fu}).  The proof uses
Malliavin Calculus. Related probabilistic formulae for spatial
derivatives of degenerate
 diffusion semigroups by Malliavin Calculus are in \cite{B} and
\cite{KS}.


 \bth \label{stime} Under Hypothesis \ref{hy},
    the following
statements hold:

 \hh (i) For any $t>0$ and   $f \in B_b ( \R^n),$ we have that $P_t f(\cdot)$
 is three times differentiable on $\R^n$ with all bounded derivatives up to
 the third order.

 \hh (ii) There exist random variables $ J_i^1 (t,x)$
  $ J^2_{ij} (t,x)$ and $J^3_{ijr} (t,x)$,
 $t>0, $ $x \in \R^n,$ $i,j, r \in \{ 1, \ldots, n \}$, which  belong
   to $L^q (\Omega)$, for any $q \ge 1$, and such that
 \begin{equation}  \label{ciao}
 \begin{aligned}
&  D_i (P_t g) (x)  =  D_i P_t g (x)   = \E [g (X_t^x)\, \, J_i^1
(t,x) ],\;\;\;
 D_{ij}^2  P_t g (x)  = \E [g (X_t^x)\, \,  J^2_{ij} (t,x) ],
\\
&  D_{ijr}^3  P_t g (x)   = \E [g (X_t^x)\, \,  J^3_{ijr} (t,x)
],\;\; \;\;\; g \in C_b (\R^n).
 \end{aligned}
 \end{equation}
(iii) For any $t >0 $, $q \ge 1,$ we have the following estimates:
\begin{equation} \label{cia}
 \begin{aligned}
&   (a) \; \E  | J_i^1 (t,x) |^q  \le  c_q (t) \, |Q_t^{-1/2}
e^{tA} e_i|^q;
\\ & (b ) \; \E  | J^2_{ij} (t,x) |^q \le c_q
(t) \, |Q_t^{-1/2} e^{tA} e_i|^q \, |Q_t^{-1/2} e^{tA} e_j|^q;
\\
& (c) \; \E  | J^3_{ijr} (t,x)|^q  \le  c_q (t) \, |Q_t^{-1/2}
e^{tA} e_i|^q \, |Q_t^{-1/2} e^{tA} e_j|^q \,|Q_t^{-1/2} e^{tA}
e_r|^q, \;\; x \in \R^n,
 \end{aligned}
 \end{equation}
 where $c_q(t)$ is a  continuous and increasing function on $[0,
\infty)$; $c_q(t)$
  $ = c (q,t, \| DF\|_{0},$ $ \| D^2 F\|_{0},$
  $ \| D^3F\|_{0},$  ${\tilde p},  \nu_1, A, n)$, where the integer
   $ \tilde p$ is introduced  in
  \eqref{q0}.
\eth \noindent It is worth noticing  that  the quantity
$|Q_t^{-1/2} e^{tA} h|^2$, corresponding to
 $q=2$,
  has a
 well known control-theoretic interpretation; see, for instance, \cite{Za}.

  Moreover,
   the following estimated are  known, see \cite{Se} and
   \cite[formula (3.4)]{Lu},
 \beq \label{st}
 |Q_t^{-1/2} e^{tA} e_i | \le \frac{ c}{ t^{h + \, 1/2} },\;\; \;\;\;
  e_i \in
 E_h
 (\R^n),\;\; 0 \le h \le k, \; t \in (0,1].
 \eeq
 where $c  \!= \! c ( {\tilde p},  \nu_1, \nu_2, A, n)>0$  and the integer
  $k$ is defined in
  \eqref {kal}.
   Estimates
 \eqref{st} can be also  deduced by purely control theoretic
 arguments. To this purpose one has to use \cite[Proposition I.1.3]{Za}
 together with \cite{AC}.
\bco \label{c3} There exists $c=  c ({\tilde p},  \nu_1, \nu_2, A,
n, \| DF\|_{0},$ $ \| D^2 F\|_{0},$
  $ \| D^3F\|_{0})>0$  such that the
following estimates hold, for any $ t>0$, $g \in B_b (\R^n),$
indices $i \in I_h $, $j \in I_{h'}$ and $r \in I_{h''}$, where
$h, h', h'' \in \{ 0, \ldots,   k\}$,
 \beq \label{chi} \bal  &&  \| D_i P_t g \|_{0}
 \, \le  c \Big ( \frac{1}{ t^{h + \, 1/2} } +  1 \Big ) \|
g\|_0; \; \;\;
  \| D_{i j}^2 P_t g \|_{0}
\, \le  c \Big ( \frac{1}{ t^{h \, + \, h' + \, 1} } + 1 \Big ) \|
g\|_0; \;\;
 \\
&& \| D_{i j r}^3 P_t g \|_{0} \, \le  c \Big (\frac{1}{ t^{h + \,
h' +\, h'' +  3/2} } + 1 \Big ) \| g\|_0. \eal \eeq \eco
 \bpf
    It is enough to prove the
 estimates when $g \in C_b (\R^n)$ (see, for instance,  \cite[Remark 3.5]{Pr}).
     Using   Theorem \ref{stime} and formula \eqref{st}, we
     first prove  the
  estimates assuming in addition that $0 < t < 1$.
We have, for any $x \in \R^n, \, t \in (0,1),$
$$
\bal | D_i P_t g (x) | \le \| g\|_0  \,  \E | J_i^1 (t,x) | \, \le
\; c_1 |Q_t^{-1/2} e^{tA} e_i | \, \| g\|_0 \le \frac {c}{ t^{h +
\, 1/2} } \| g\|_0.
 \eal $$
 In a similar way, we  get   the second and third  estimates, for
 $t < 1$.

 When  $t \ge 1$,   by the semigroup and
 the
 contraction property of $P_t$, we have:
$$
\| D_i P_t g \|_{0}  = \| D_i P_{1/2} (P_{ t - \frac{1}{2} } \, g)
\|_{0} \le c 2^{h + 1/2} \|  P_{\frac{ 2t - 1 }{2}} \,  g \|_{0}
\le  c 2^{h + 1/2}
  \|  g  \|_{0},
$$
for any $0 \le i \le k$. Hence the required estimate  of $D_i P_t
g$ follows for any $t>0$.  Similarly, we get the other estimates
for any $t>0$. \epf
 The main result of the section is the following theorem.
  Its
  proof
  also allows
  to   complete the final part of the proof  of \cite[Theorem 3.4]{Lu}.
 We set $t \wedge 1 = \min (t,1)$.
 \bth \label{sti}
 Fix any $ \gamma \in (0,3) $ {\it non-integer}. There exists $c=
    c ( \gamma, {\tilde p},  \nu_1, \nu_2, A,
n, \| DF\|_{0},$ $ \| D^2 F\|_{0},$ $ \| D^3F\|_{0})>0$, such
that, for any
  $f \in {\cal C}^{\gamma}_d
  (\R^n)$, $t>0$,  for any  indices  $i
\in I_h $, $j \in I_{h'}$ and $r \in I_{h''}$, where  $h, h', h''
\in \{ 0, \ldots,   k\}$, it holds
 \begin{align} \label{c7}
  \nonumber & (i)\;
  \| D_i P_t f \|_{0} \le c \Big ( \frac{1}{  (t\wedge 1 )^{  \frac{1 - \gamma}{2} + h
 } }+ 1 \Big ) \| f\|_{\gamma, d};
  \;\;\;     (ii)\,
 \,  \| D_{ij}^2 P_t f \|_0
\le c \Big ( \frac{1}{ (t\wedge 1 )^{ h  + h' +
 \frac{ \, 2 - \gamma}{2}
} } + 1 \Big ) \| f\|_{\gamma, d};
\\
   &  (iii)\;
 \;  \| D_{ij r }^3
P_t f \|_0 \le c \Big (\frac{1}{t^ { h \, + h' + \, h'' + \,
\frac{3 - \gamma}{2} }  } + 1 \Big )\| f\|_{\gamma, d};
 \;\;\;\;\;
  (iv) \;\; \| P_t f \|_{\gamma, d} \le c \| f\|_{\gamma, d}.
 \end{align}
 \eth

 \begin{remark} \label{elli} { \em Estimates (i)-(iv) will be used
 to get elliptic and parabolic Schauder estimates for $\cal A$.
  However, we
 stress that
 to prove  {\it  elliptic Schauder estimates} we only  need a
 special
 case of \eqref{c7}.
 More precisely, we need, for any $\theta \in (0,1)$,
 $f \in {\cal C}^{\theta}_d (\R^n),$ $t>0$,
  for any indices $i,\, j \in I_0, \, r \in I_h$,
 with $h \in \{ 0, \ldots,   k\}$,
 \beq \label{elli1}
\bal  & (a)\; \| D_r P_t f \|_{0} \le c \Big ( \frac{1}{  t^{
\frac{1 - \theta}{2} + h
 } }+ 1 \Big ) \| f\|_{\theta, d};
 \;\; \;
 (b)\; \| D_{ij}^2 P_t f \|_0 \le c \Big
( \frac{1}{ t^{ \frac{2 - \theta}{2}  } } + 1 \Big ) \|
f\|_{\theta, d};
\\
& (c) \; \| D_{ij r }^3 P_t f \|_0 \le c \Big (\frac{1}{t^ {
\frac{3 - \theta}{2} + h }  } + 1 \Big )\| f\|_{\theta, d};
 \;\;\; (d) \; \| P_t f \|_{\theta, d}
 \le c \| f\|_{\theta, d}.
  \eal \eeq
 These estimates are  simpler to obtain than the general
  ones  in which $\gamma \in (0,3)$.
 On the other hand, the estimates (iv) in \eqref{c7}
 with $\gamma \in (2,3)$
 are a particular case of
    parabolic Schauder estimates corresponding to $H=0$ in
     \eqref{2} (see Theorem \ref{buon1}). Estimates (iv)
      will be deduced by (iii).
 \qed} \end{remark}
  In order to prove the main result we need three  preliminary
  lemmas. To state the first one
   we introduce the deterministic process $Y_t^x $,
 \beq \label{go}
 Y_t^x = e^{tA} x + \int_0^t e^{(t-s) A} F(Y_s^x)ds,\;\; t \ge 0, \; x \in
 \R^n,
 \;\; \mbox{ which solves} \;\;
\left\{
\bal \dot  {Y_t } & = AY_t + F(Y_t), \\
      Y_0 & = x,
\eal \right. \eeq
 \ble \label{g}
 For any $q>0$,  there exists  $C = C(q, {\tilde p},
 \nu_1, \nu_2, n, A,  \| D F \|_0)>0$, such that
 \beq \label{che1}
 \sup_{x \in \R^n} \, \E [ \big( d (  X_t^x  , Y_t^x )
\big) ^{q}]
 = \sup_{x \in \R^n}  \E
{\nor}  X_t^x  - Y_t^x   {\nor}^q  \, \le \, C \, t^{\frac{q}{2}
}, \;\;\;\;\;\;\;\; 0 \le t \le 1.
 \eeq
 \ele
\bpf Note that \eqref{che1} is equivalent to the following
 assertion: for any $q>0$,
$0 \le h \le k$,
  there exists  $C_1>0$   such that
 \beq \label{che}
 \sup_{x \in \R^n} \,  \E \big | E_h (X_t^x  - Y_t^x ) \big|^q  \,
\le \, C_1 \, t^{\frac{q}{2} (2h+1) }, \;\;\; 0 \le t \le 1,
 \eeq
 see \eqref{di}. Let us prove \eqref{che}.
  Since  there exists $c>0$, such
that $|x| \le c \sum_{h =0}^k |E_h x|$, for any $x \in \R^n$, we
get
 \beqr
 | E_h ( X_t^x - Y_t^x) | \le \Big |   \int_0^t  (E_h e^{(t-r) A} E_0) \,    [ F(X_r^x)- F
 (Y_r^x)]dr \Big | + | E_h Z_t^0|
\\ \le     c  \| DF \|_0 \,   \sum_{j =0 }^k \int_0^t  \| E_h
e^{(t-r) A} E_0 \|_L \, | E_j (X_r^x
 - Y_r^x ) | dr \,  + \, | E_h Z_t^0|,
 \eeqr
  $\P$-a.s.. Using   the  following estimate, see \cite[Lemma
 3.1]{Lu},
 \beq \label{l}
 \| E_h e^{s  A} E_0
\|_L \le c'' s ^{h },\;\; s \in [0,1],\;\; 0 \le h \le k,
 \;\;\; \mbox{where}\; c'' = c''(A)>0,
 \eeq
we arrive at
 \begin{align} \label{rit}
\nonumber | E_h ( X_t^x - Y_t^x) | \le | E_h Z_t^0|
 \\  + \,   C  \int_0^t  (t -
 r)^h  \, | E_h (X_r^x
 - Y_r^x) | dr \,
+  C  \sum_{j=0,\;  j \not =  h }^k \int_0^t   (t -
 r)^h  \, | E_j (X_r^x
 - Y_r^x) | dr.
 \end{align}
 $\P-\text{a.s.}.$ Now we use  that $| E_j (X_r^x
 - Y_r^x) | \le |X_r^x
 - Y_r^x |$, $0\le j \le k$. Since
  $$
 |X_t^x  - Y_t^x  | \le C' \int_0^t |X_s^x  - Y_s^x  | ds +
 |Z_t^0|,
$$
an application of the
 Gronwall lemma gives, $\P$-a.s.,
 \begin{equation} \label{rito}
|X_t^x  - Y_t^x  | \le |Z_t^0| + C' \int_0^t |Z_s^0| e^{(t-s) C'}
ds \le |Z_t^0| + C_1 \int_0^t |Z_s^0| ds ,\;\; 0 \le t \le 1.
\end{equation}
Using estimate \eqref{rito} in \eqref{rit} we get
 \beqr
&& | E_h ( X_t^x - Y_t^x) | \\ && \le | E_h Z_t^0|
   \, +  \,   C   \int_0^t  (t -
 r)^h  \, | E_h ( X_r^x - Y_r^x)|  dr \,
+  C    \sum_{ j \not =  h } \int_0^t   (t -
 r)^h  \, \Big [ |Z_r^0| + C_1 \int_0^r |Z_s^0| ds     \Big ] \,  dr,
\eeqr $\P$-a.s.. Let now $q \in \Z_+ $  and recall that $0 \le t
\le 1.$ We have
 \beqr
| E_h ( X_t^x - Y_t^x) |^q  \\ \le C_3 \Big ( | E_h Z_t^0|^q
    +        \int_0^t    | E_h ( X_r^x - Y_r^x)|^q  dr \,
+      \sum_{ j \not =  h } \int_0^t   (t -
 r)^{hq}  \, \Big [ |Z_r^0|^q +  \int_0^r |Z_s^0|^q ds     \Big ]
 dr \Big ),
\eeqr
 $\P$-a.s.. Before applying
   the expectation in the last formula,
  we   check that
  \beq \label{oo}
 \E \big | E_h Z^0_t \big |^{{ q} } = \E  \Big |    \int_0^t E_h e^{(t-s)
 A}  Q^{1/2}
 dW_s \Big |^{q }  \le \, c_{q, h} \, t^{q (2h+1) /2},\;\; \;\;
  0\le  t \le 1,\; q>0, \; 0 \le h \le
 k,
 \eeq
 where $E_h$ are the orthogonal projections introduced in \eqref{sp}.
  Denoting by $N(0,Q_t) $  the Gaussian measure on $\R^n$ with mean
  $0$
and covariance matrix $Q_t$ given in \eqref{qtt}, we have:
 \beq \label{ri}
 \bal & \E \big | E_h Z^0_t \big |^{{ q} } =  \int_{\R^n}
 |E_h y|^{q } \, N(0, Q_t) dy
 \\ &   = \int_{\R^n}
 |E_h Q_t^{1/2} z|^{q } \, N(0, I) dz
 \le
  \| E_h Q_t^{1/2} \|_L^q  \, \int_{\R^n}
  | z|^{q } \, N(0, I) dz
\le  c \, t^{q(2 h +1 )/2}, \;\; t \le 1,
 \eal \eeq
where $I$ is the $n \times n$ identity matrix. In the last
inequality we have used that $\| E_h Q_t^{1/2} \|_{ L}   \le $ $c'
t^{(2h +1)/2 }$, $0 \le t \le 1$, $0 \le h \le k$,
 where $c' = c'({\tilde p},n,A, \nu_1, \nu_2)$
 (see \cite[ formula
(3.2)]{Lu}).

 By \eqref{oo}, we infer
 \beqr \E | E_h ( X_t^x - Y_t^x) |^q
\le C_3 \Big ( \E | E_h Z_t^0|^q
    +        \int_0^t     \E | E_h ( X_r^x - Y_r^x)|^q  dr \,
\\ +      \sum_{ j \not =  h } \int_0^t   (t -
 r)^{hq}  \, \Big [ \E |Z_r^0|^q +  \int_0^r \E |Z_s^0|^q ds     \Big ]
 dr \Big )
\\
\le C_4 \Big ( t^{\frac{q(2h+1)}{2}}
    +        \int_0^t     \E | E_h ( X_r^x - Y_r^x)|^q  dr \,
 +      \sum_{ j \not =  h } \int_0^t   (t -
 r)^{hq}  \, \big [ r^{q/2} +  \frac{r^{1 + q/2}} {1 + q/2}   \big ]
 dr \Big ).
\eeqr
 Using that $\int_0^t (t-s)^p s^{r} ds = \frac{p!}{ (r+p+1) (r+p) \ldots
  (r+1)} t^{r+p+1}$, for $p \in \Z_+$, $r>0$,  we get
\beqr
 \E | E_h ( X_t^x - Y_t^x) |^q
\le C_5 \Big ( t^{\frac{q(2h+1)}{2}}
    +        \int_0^t     \E | E_h ( X_r^x - Y_r^x)|^q  dr \,
 +    2   \sum_{ j \not =  h } \int_0^t   (t -
 r)^{hq}  \,  r^{q/2}
 dr \Big )
\\
 \le C_6 \Big ( t^{\frac{q(2h+1)}{2}}
    +        \int_0^t     \E | E_h ( X_r^x - Y_r^x)|^q  dr \,
 +     t ^{hq + 1+ q/2}
 \Big ) \\ \le  2 C_6 \Big ( t^{\frac{q(2h+1)}{2}}
    +        \int_0^t     \E | E_h ( X_r^x - Y_r^x)|^q  dr \,
 \Big ), \;\; t \le 1.
\eeqr
 Applying the Gronwall lemma, we get
  $$
 \E | E_h ( X_t^x - Y_t^x) |^q  \le C_7 \, \, t^{\frac{q(2h+1)}{2}},\;\;
 0 \le t \le 1.
 $$ Now if $q \in \R_+$, $q>0$, we consider an integer $m \ge q$.
 By the Jensen inequality,
$$
 \Big (\E | E_h ( X_t^x - Y_t^x) |^q \Big )^{m/q}
 \le \, \E | E_h ( X_t^x - Y_t^x) |^m \le
 C t^{\frac{m(2h+1)}{2}}, \;\; t \le 1.
$$ This implies that $\E | E_h ( X_t^x - Y_t^x) |^q \le C ^{q/m}
 \, t^{\frac{q(2h+1)}{2}}$.
 The assertion is proved.
\epf
 \ble \label{ci}  For any $\omega$, $\P$-a.s., $t \in [0,1]$, the mapping
  $x \mapsto X_t^x (\omega) \in \R^n$ is differentiable up to the third order
  on $\R^n$.
   Moreover, for any $i, j , r \in \{ 0 , \ldots, n\}$, $x \in \R^n$,
     there
   exist
     continuous adapted stochastic processes $(\eta_i (t,x))$,
  $(\eta_{ij}(t,x))$ and  $(\eta_{ijr}(t,x))$ with values in
  $\R^n$ and $C=  C ( \| D F\|_0, \| D^2 F\|_0, \| D^3 F\|_0, \| A\|_{L}\,  )
> 0$   such that
 \beqr
 \eta_i (t,x) = D_i \, X_t^{x} = \lim_{h \to 0} \, { (X_t^{x + h
e_i} - X_t^x)} {\, h^{-1}},\;\;\;\;\;\; \eta_{ij} (t,x) = D_{ij}^2
\, X_t^{x},\;\;\;\;\;\; \eta_{ijr} (t,x) = D_{ijr}^3 \, X_t^{x}
\\
\mbox{and} \;\;\;\; |\eta_i (t,x)| + | \eta_{ij}(t,x) | +  |
\eta_{ijr}(t,x)|
 \le C,\;\;\; \mbox { for any  $ t \in [0,1],\; x \in \R^n,
  \; \omega \in \Omega,
 \; \P$-a.s.}.
 \eeqr
 \ele
 \bpf
  The proof is  straightforward. We include it for the sake
  of completeness.
   Fix $\omega \in
 \Omega$, $\P$-a.s., and
 introduce the Banach space $ E = C([0,1]; \R^n)$.  Define  the map
${\cal F} : \R^n \times E \to E$,
$$
{\cal  F} (x, u)(t)  :=  u(t) - x - \int_0^t \big ( A u(r) \, + \,
F(u(r)) \big )\, d r - \sqrt{Q} W_t(\omega) , \;\;\; t \in [0,1],
\; u \in E,\; x \in \R^n.
$$
Applying  the implicit function theorem, we find that  the
mapping:
  $x \mapsto X^x_{(\cdot)} (\omega)$ from $\R^n$  into $E$ is three times
    Fr\'echet-differentiable.
 Denote by $\eta_i (t,x)$, $\eta_{ij} (t,x) $ and $\eta_{ijr}
 (t,x)$
  $t \in [0,1]$, respectively  the
 first (directional)
 derivative at $x \in \R^n$ in the direction $e_i$,  the second derivative
 at $x$ in the directions $e_i$ and $e_j$, and the third derivative
 at $x$ in the directions $e_i$, $e_j$ and $e_r$,
 where $i, j, r =1, \ldots, n$. Note that $\eta_i (t,x)$, $\eta_{ij} (t,x) $ and $\eta_{ijr}
 (t,x)$ solves, $\P$-a.s.,
  the variation
   equations
  \beqr
&&  \eta_i (t,x) = e_i +  \int_0^t \big ( A \eta_i (s,x) \, + \, D
F (X_s^x) [\eta_i (s,x)] \big) ds;
\\
&& \eta_{ij} (t,x) =  \int_0^t \big ( A \eta_{ij} (s,x) \, + \,
D^2 F (X_s^x) [\eta_i (s,x)] \, [\eta_j (s,x)] \, + \, D F (X_s^x)
[\eta_{ij} (s,x)] \big) ds;
\\
&& \eta_{ij r} (t,x) =  \int_0^t \big ( A \eta_{ijr} (s,x) \, + \,
D^3 F (X_s^x) [\eta_i (s,x)] \, [\eta_j (s,x)] \, [\eta_r (s,x)]
\big) ds \,
\\
&&   + \int_0^t \big ( D^2 F (X_s^x) [\eta_{ir} (s,x)] \, [\eta_j
(s,x)]  + D^2 F (X_s^x) [\eta_i (s,x)] \, [\eta_{jr} (s,x)]
 \, + \, D F (X_s^x) [\eta_{ijr} (s,x)] \big) ds,
 \eeqr
 $ t\in [0,1].$
  It follows easily that   $(\eta_i (\cdot ,x))$,
  $(\eta_{ij}(\cdot,x))$ and  $(\eta_{ijr}(\cdot,x))$ are
   continuous adapted stochastic processes.  An application of the
   Gronwall lemma gives the final assertion.
 \epf
 \ble \label{tec} Let $f
\in {\cal C}^{\gamma}_d (\R^n)$, $\gamma \in (2,3)$, and $i,  j, r
\in \{ 1, \ldots, n\}$. Consider the following random variables
depending on $t \in (0,1)$ and $x \in \R^n$
 (see
 (\ref{d0}) and \eqref{ciao}$)$ fff
$$
\Lambda (t,x) = \lan D_{E_0 \,} f (Y_t^x), E_0 (X_t^x - Y_t^x)
\ran
 +  \frac{1}{2} \lan D^2_{E_0 \,}  f (Y_t^x) \, [E_0 (X_t^x - Y_t^x)],
  E_0 (X_t^x - Y_t^x) \ran.
$$
Then the  functions:
   $ \ds { \phi_i (x,t) = \E \Big [ \Lambda (t,x) \, J^1_{i} (t,x)
 \Big ],}$
 $  \phi_{ij} (x,t) = \E \Big [ \Lambda (t,x) \,  J^2_{ij} (t,x)
  \Big ], $ $
  \phi_{ijr} (x,t) $ $= \E \Big [ \Lambda (t,x) \,  J^3_{ijr} (t,x)
 \Big ],
 $
 $x \in \R^n, \, t \in (0,1)$,
 are continuous and bounded on $\R^n \times (0,1)$.\ele

 \bpf

\vskip 1mm Let us  treat $\phi_i$. We introduce the deterministic
functions $K :
  \R^n \times \R^n \times [0,1] \to \R $,
 \beq \label{gg}
K (x,z,t) =   \lan D_{E_0 \,} f (Y_t^z), E_0 ( x \, - \, Y_t^z)
\ran
 +  \frac{1}{2} \lan D^2_{E_0 \,}  f (Y_t^z) \, [E_0 ( x \, - \, Y_t^z)],
  E_0 (x \, - \, Y_t^z) \ran
 \eeq
and $g_i:  \R^n \times \R^n \times (0,1] \to \R$,
 $$
g_i (x,z,t) = \E \Big [ K ( X_t^{x}, z, t) \, J^1_{i} (t,x)
 \Big ],\;\;\; x,z \in \R^n,\; t \in (0,1].
$$
 Note that $\phi_i (x,t) = g_i (x,x,t)$, $x \in \R^n$, $t \in
 (0,1)$. We first prove that
 \beq \label{tec1}
 g_i(x,z, t) =  D_i \big (  \E \big [ K ( X_t^{(\cdot)}, z, t) \big ] \big
 )(x)= \E \big [
  \lan D_{x}K ( X_t^{x}, z, t),
   \eta_i (t,x) \ran \big ],
 \eeq
$x, z \in \R^n, \; t \in (0,1)$ (here $D_{x_i} = D_i$  denotes the
partial derivative with respect to $e_i$
 and $D_x$  denotes the gradient in the $x$-variable;
$\eta_i$ is introduced in Lemma \ref{ci}). To this purpose, remark
that it holds
 \beq \label{ft}
|K(x,z,t)| +   | D_{x_i} K(x,z,t) | + | D_{x_i x_j}^2 K(x,z,t) | +
| D_{x_i x_j x_r}^3 K(x,z,t) | \le 8 \| f\|_{\gamma, d} \, ( 1 +
|E_0 (x  -  Y_t^z)|^2),
 \eeq
$t \in [0,1],\; x , \, z \in \R^n,$ $i,j,r \in \{1, \ldots, n\}$.
Moreover, an
 application of the Gronwall lemma shows that
\beq \label{f7}
  |X_t^x|  \, \le \, e^{(\| A\|_{L} \, + \| DF\|_0)}
   \, \, \big (|x| \, + \, |F(0)| \, + \,
     \| \sqrt{Q}\|_{L} \, \sup_{s \le 1}|W_s|  \big
 ),\;\;\;  t \in [0,1],\; x  \in \R^n,
\eeq
 $\P$-a.s..
 By  \eqref{ft} and \eqref{f7}, using   Lemma \ref{ci}, we
  get the existence of the   partial derivatives
 $$
 D_{x_i} \big (  \E \big [ K ( X_t^{(\cdot)}, z, t) \big ] \big
 )(x) =   \E \big [ \lan D_{x} K ( X_t^{x}, z, t), \, \eta_i (t,x)
  \ran \big ],
 \;\; x, z \in \R^n, \; t \in (0,1),\;  1 \le i \le n.
$$ To
 obtain \eqref{tec1}, we consider  test functions $\varphi_m \in C_0^{\infty}
  (\R^n)$ such that $0 \le \varphi_m \le 1$, $m \in \N $, $\varphi_m (x) =1
  $, when $|x| \le n$, $\varphi_m (x) =0
  $, when $|x| > m+1$ and $|D \varphi_m (x)| \le 1$, for $x \in
  \R^n$, $m \in \N.$  By Theorem \ref{stime} and Lemma \ref{ci}, we know that, for $
x, z \in \R^n, \; t \in (0,1),$ $m \in \N$,
  \beqr
  D_{x_i} \big (  \E \big [ K ( X_t^{(\cdot)}, z, t) \, \varphi_m
(X_t^{(\cdot)})\big ] \big
 )(x)  =  \E \big [ K ( X_t^{x}, z, t) \, \varphi_m (X_t^{x}) \,
  J^1_i (t,x) \big ]
  \\  = \E \big [
   \lan D_{x }K ( X_t^{x}, z, t), \eta_i (t,x) \ran
    \varphi_m (X_t^x) \, + \,
  K ( X_t^{x}, z, t) \ \lan D \varphi_m (X_t^x),
   \, \eta_i (t,x) \ran \big ]
\eeqr Passing to the
 limit as $m \to \infty$, we get \eqref{tec1}, by the dominated convergence
 theorem.
  By \eqref{tec1},  \eqref{ft} and  Lemma \ref{ci}, we have (setting $z=x$)
   \beqr
 |\phi_i(x,t)| = | g_i(x,x, t) | \le \E \big [ |
  D_{x_i}K ( X_t^{x}, x, t)|
  \, |\eta_i (t,x)| \big ]
 \\ \le 8 \| \eta_i\|_{L^{\infty}} \,  \| f\|_{\gamma, d} \, ( 1 +
\E [|E_0 (X_t^x \, - \, Y_t^x)|^2]),
 \eeqr for any $ x \in \R^n$,
$t \in (0,1)$. Applying   Lemma \ref{g}, we get
$$
\sup_{x \in \R^n,\, t \in (0,1)} \, |\phi_i (x,t)|
 \, \le \,  8  c \|
\eta_i\|_{L^{\infty}} \, \| f\|_{\gamma, d} \, \sup_{ t \in (0,1)}
| 1 + t| \le 16 c  \, \| \eta_i\|_{L^{\infty}} \, \| f\|_{\gamma,
d}.
 $$
 To treat $\phi_{ij}$ and $\phi_{ijr}$ we proceed similarly.
  Concerning  $\phi_{ij}$ we  introduce
  \beqr
g_{ij} (x,z,t)=  \E \big [ K ( X_t^{x}, z, t) \, J^2_{ij} (t,x)
\big] = D_{ij}^2 \big (  \E \big [ K ( X_t^{(\cdot)}, z, t) \,\big
] \big
 )(x)
\\
= \E \big [
  \lan D_{x}^2K ( X_t^{x}, z, t)
  \, [\eta_j (t,x)],  \eta_i (t,x) \ran \, + \,
  \lan D_{x}K ( X_t^{x}, z, t),
  \, \eta_{ij} (t,x) \ran \big ].
   \eeqr
Since $\phi_{ij}(x,t) = g_{ij}(x,x,t)$, we obtain the assertion
for $\phi_{ij}$, using  \eqref{ft}, Lemmas \ref{ci} and \ref{g},
as before.  To treat $\phi_{ijr}$  we introduce
 $g_{ijr} (x,z,t)=  \E \big [ K ( X_t^{x}, z, t) \, J^3_{ijr} (t,x)
 \big]$. Note that
 \beqr
g_{ijr} (x,z,t) \!=\!  \E \big [
  \lan D_{x}^3 K ( X_t^{x}, z, t)
   [\eta_r (t,x)] [\eta_j (t,x)],  \eta_i (t,x) \ran  \! + \!
  \lan D_{x }^2 K ( X_t^{x}, z, t) [ \eta_{jr} (t,x)],  \eta_i
  (t,x) \ran
\\ + \lan D_{x }^2 K ( X_t^{x}, z, t) [\eta_j (t,x)],  \eta_{ir}
(t,x) \ran    +  \lan D_{x }^2 K ( X_t^{x}, z, t)
   [ \eta_{r} (t,x)],  \eta_{ij} (t,x) \ran
   \\  +
    \lan D_{x}K ( X_t^{x}, z, t)
  , \eta_{ijr} (t,x) \ran \big ].
 \eeqr
Since   $\phi_{ijr} (x,t)= g_{ijr}(x,x,t)$, we get the assertion
 for $\phi_{ijr} $ proceeding as for $\phi_i$ and $\phi_{ij}$.
 The proof is complete. \epf

\hh {\bf Proof of Theorem \ref{sti}.}  \
 Thanks to Corollary \ref{c3}, {\it it is enough to prove all the
estimates for $ 0 < t < 1. $
 }
 Indeed,  concerning \eqref{c7}, we have, for $t \ge 1$,
 \beqr   \| P_t f \|_{\gamma , d} =
  \sum_{m = 0}^k \;  \sup_{z \in \R^n} \, \|
 (P_{t} \,
f) (z + \cdot) \|_{C^{\gamma / (2m +1)}_b \, (E_m (\R^n))} \,
\\ \le \,  c'  \, \sum_{m = 0}^k \, \big ( \| f\|_0 \, + \,
\sup_{z \in \R^n} \, \|
 (P_{t} \,
f) (z + \cdot) \|_{C^{3}_b \, (E_m (\R^n))} \, \big )  \,  \le \,
c_2 \| f \|_0,\;\;\; \;\;  f \in {\cal C}^{\gamma}_d(\R^n).
 \eeqr
 {\it We  will
  show the estimates only for $\gamma \in (2,3)$
 non-integer.}

\vskip 1mm \noindent Indeed, the  cases of $\gamma \in (0,1)$
 and $\gamma \in (1,2) $  can be similarly  treated
  and are even simpler.
  Alternatively, once we have proved the
  estimates  for $\gamma \in (2,3)$, the
   remaining estimates can be obtained  by an  interpolation argument.
  Let us briefly explain such method  which has been also used
   in the proof of  \cite[Theorem 3.4]{Lu}.
 We assume that (i)-(iv) hold for $\gamma = 5/2$ and show that
 they hold also for a fixed $\gamma' \in  (0,2)$ non-integer.
 By \cite[Theorem 2.2]{Lu},
  we know in particular that
 \beq \label{f1}
 \big ( C_b (\R^n), {\cal C}^{5/2}_d (\R^n) \big)_{2{\gamma '} / 5,\, \infty}
 \, =\, {\cal C}^{{\gamma '}}_d (\R^n).
 \eeq
To be precise, \eqref{f1} is proved in \cite{Lu} when $C_b (\R^n)$
denotes the Banach space of all real continuous and bounded
functions defined on $\R^n$. However,  the same proof of \cite{Lu}
works also when we consider $ C_b (\R^n)$ as the space of
 all real uniformly  continuous and bounded functions.
  Concerning estimate (iv) in
   \eqref{c7}, by \eqref{f1} and  \cite[Proposition 1.2.6]{Lu2} we get
$$
 \| P_t \|_{L ( {\cal C}^{{\gamma '}}_d (\R^n), {\cal C}^{{\gamma '}}_d (\R^n)  ) } \,
 \le \,  (\| P_t \|_{L ( C_b (\R^n), C_b (\R^n)  ) })^{1 - \frac{2 {\gamma '}} {5} } \,
 (\| P_t \|_{L ( {\cal C}^{5/2}_d (\R^n), {\cal C}^{5/2}_d (\R^n)  ) })^{\frac{2
 {\gamma '}}
 {5}} \le C,
$$
for $t \ge 0$.  As for (iii), we fix $x \in \R^n$, $t \in (0,1]$
and introduce the linear
 operator $T_{x,t} : {\cal C}^{{\gamma '}}_d (\R^n) \to \R$, $T_{x,t} f : =
  D_{ijr}^3 P_t f (x) $, for any $f \in {\cal C}^{{\gamma '}}_d(\R^n)$. We
  have:
$$ \| T_{x,t} \|_{L ( {\cal C}^{{\gamma '}}_d (\R^n), \R )} \,
 \le \,  (\| T_{x,t} \|_{L ( C_b (\R^n), \R  ) })^{1 - \frac{2 {\gamma '}} {5} } \,
 (\| T_{x,t} \|_{L ( C^{5/2}_d (\R^n), \R  ) })^{\frac{2 {\gamma '}}
 {5}} \le c t^{-( \frac{3}{2} + h + h' + h'') + \frac{\gamma ' }{2}
 },
$$ $t \in (0,1]$ (uniformly in  $x\in \R^n$). In  a similar way, one can prove (i)
and (ii) for $\gamma'$.

 \hh {\it - \  We prove the first  estimate in \eqref{c7}, for $t
\in (0,1)$, $\gamma \in (2,3)$ non-integer and $i \in I_h$.}

\hh We start from \eqref{ciao} and write
 \beq \label{fg}
 \bal
 D_i P_t f (x)   = \Lambda_1(t,x) + \Lambda_2 (t,x) \;\; \mbox{where}
\\
 \mbox{ $\Lambda_1(t,x) =  \E \big [ \{ f(X_t^x) -   f \big ( E_0 X_t^x +
 \sum_{h=1}^k E_h Y_t^x \big ) \} \,  J^1_{i} (t,x) \big ]; \;\;$ }
\\  \mbox{ $ \Lambda_2 (t,x) = \E \big [ f \big ( E_0 X_t^x +
 \sum_{h=1}^k E_h Y_t^x \big ) \,  J^1_{i} (t,x)
 \big ],$}
\eal \eeq
 where $(Y_t^x)$ is defined in \eqref{go}. Let us treat $\Lambda_1$ and  $\Lambda_2$
 separately. We have since $0 < \gamma / (2 m +1 ) < 1$ if $m = 1 , \ldots, k$
 (using \eqref{cia}, \eqref{st} and Lemma \ref{g})
\beqr
  |\Lambda_1(t,x)| \le c \| f\|_{\gamma, d} \, \E \Big [ \big \{
\mbox{ $ \sum_{m=1}^k $}  |E_m (X_t^x - Y_t^x)|^{\frac{\gamma}{2m
+1}} \big \} \,\, |J^1_{i} (t,x)| \Big]
\\
  \le c' \| f\|_{\gamma, d} \Big(\E \, { \nor X_t^x - Y_t^x
\nor }^{ {2 \gamma}}  \Big )^{1/2}  \,
 \big (\E |J^1_{i} (t,x)|^2 \big)^{1/2}
\le c_2 \| f\|_{\gamma, d} \, t^{\frac{ \gamma}{2 } } t^{-(h+  1/2
)} = c_2 \| f\|_{\gamma, d} \,\, t^{\frac{ \gamma -1}{2 } -h },
 \eeqr
 $t \in (0,1)$, uniformly in $x \in \R^n$. Let us concentrate on the more
difficult term $\Lambda_2$.
 We  write
 \beq \label{lambda2}
 \bal  \Lambda_2 (t,x)  & = \Lambda_{21} (t,x) + \Lambda_{22} (t,x),
 \;\; \mbox{where}
\\
 \Lambda_{21} (t,x) & =  \E \Big [  \Big ( f \big ( E_0 X_t^x +
 \mbox{$ \sum_{m=1}^k $} E_m Y_t^x \big ) - f(Y_t^x) -
 \lan D_{E_0 \,} f (Y_t^x), E_0 (X_t^x - Y_t^x) \ran
\\ & - \mbox{$ \frac{1}{2}$} \lan D^2_{E_0 \,}  f (Y_t^x) \, [E_0 (X_t^x - Y_t^x)],
  E_0 (X_t^x - Y_t^x) \ran
\Big )\,  J^1_{i} (t,x)
 \Big ],
\\
 \Lambda_{22} (t,x) & =
 \E \Big [ \Big ( \lan D_{E_0 \,} f (Y_t^x), E_0 (X_t^x - Y_t^x) \ran
 \\ & +  \mbox{$ \frac{1}{2} $} \lan D^2_{E_0 \,}  f (Y_t^x) \, [E_0 (X_t^x - Y_t^x)],
  E_0 (X_t^x - Y_t^x) \ran
\Big )\,  J^1_{i} (t,x)
 \Big ],
 \eal
 \eeq
see \eqref{d0}. Note that, since $(Y_t^x)$ is deterministic, $\E [
f(Y_t^x)J^1_{i} (t,x) ]= f(Y_t^x) \E [ J^1_{i} (t,x) ] = f(Y_t^x)
D_i (P_t 1) (x) =0$, for any $x \in \R^n$, $t >0$, $1 \le i \le
n$.

\vskip 1mm To estimate $\Lambda_{21}$, remark that
 $f(x+  \cdot) \in
 C^{\gamma}_b (E_0 (\R^n))$, $\gamma \in (2,3)$, uniformly in $x$.
 By the
mean value theorem, we have:
 \beqr
&& \sup_{x \in \R^n}\, |\Lambda_{21} (t,x)| \le \| f\|_{\gamma, d}
\, \sup_{x \in \R^n}  \, \{ \E \big [ | E_0 (X_t^x - Y_t^x ) |^{ {
\gamma}} \, |J^1_{i} (t,x) |\big ] \, \}
\\ && \le
\| f\|_{\gamma, d} \sup_{x \in \R^n}  \, \Big( \E \big [ | E_0
(X_t^x - Y_t^x ) |^{ {2 \gamma}} \Big ] \Big )^{1/2} \, \cdot \,
\sup_{x \in \R^n}  \, \big (\E |J^1_{i} (t,x)|^2 \big)^{1/2} \le
c_3 \| f\|_{\gamma, d} \; t^{\frac{ \gamma -1}{2 } -h },
 \eeqr
see also \eqref{che}. Finally, using Lemma \ref{tec}, we infer $
\sup_{x \in \R^n, \; t \in (0,1)} \; | \Lambda_{22} (t,x)| = c_4 <
\infty. $ This proves the  estimate.

\hh {\it - \  We prove   (ii) and  (iii) in \eqref{c7}, for $t \in
(0,1)$ and $\gamma \in (2,3)$ non-integer.}

\hh These estimates can be similarly proved  to the first estimate
 in \eqref{c7}.
 We only give
 the proof of (ii). Let $i \in I_h$ and $j \in I_{h'}$.
  We  write
 $$
 \bal
& D_{ij}^2 P_t f (x)   = {\tilde \Lambda}_1(t,x) + {\tilde
\Lambda}_{21} (t,x) \,
 + \, {{\tilde \Lambda}}_{22} (t,x),
 \;\; \mbox{where}
\\ \\
& \mbox{$ {\tilde \Lambda}_1(t,x) =  \E \big [ \{ f(X_t^x) -   f
\big ( E_0 X_t^x +
 \sum_{h=1}^k E_h Y_t^x \big ) \} \,  J^2_{ij} (t,x) \big ], \;\;
 $}
\\ & {\tilde \Lambda}_{21} (t,x) =  \E \Big [ \Big ( f \mbox { $
\big ( E_0 X_t^x +
 \sum_{h=1}^k E_h Y_t^x \big ) $}
- f(Y_t^x) - \lan D_{E_0 \,} f (Y_t^x),
 E_0 (X_t^x - Y_t^x) \ran
\\  & - \frac{1}{2} \lan D^2_{E_0 \,}  f (Y_t^x) \, [E_0 (X_t^x - Y_t^x)],
  E_0 (X_t^x - Y_t^x) \ran
\Big )\,  J^2_{ij} (t,x)
 \Big ],
\\
& {\tilde \Lambda}_{22} (t,x)  =
 \E \Big [ \Big ( \lan D_{E_0 \,} f (Y_t^x), E_0 (X_t^x - Y_t^x) \ran
\\ & +  \frac{1}{2} \lan D^2_{E_0 \,}  f (Y_t^x) \, [E_0 (X_t^x -
Y_t^x)],  E_0 (X_t^x - Y_t^x) \ran \Big )\,  J^2_{ij} (t,x)
 \Big ],
\eal $$
 $t \in (0,1), \, x \in \R^n$.
  We have (using \eqref{cia}, \eqref{st},  Lemmas \ref{g} and \ref{tec})
\beqr
  \sup_{x \in \R^n} \, |{\tilde \Lambda}_1(t,x)| \le c \| f\|_{\gamma, d} \,
 \sup_{x \in \R^n}  \,  \E \Big [ \big \{
\sum_{m=1}^k  |E_m (X_t^x - Y_t^x)|^{\frac{\gamma}{2m +1}} \big \}
\,\, |J^2_{ij} (t,x)| \Big] \le  c_2 \| f\|_{\gamma, d} \,
t^{\frac{ \gamma - 2} {2} \, - \, h \, - \, h' }.
 \eeqr
By the mean value theorem, we find
 \beqr
 \sup_{x \in \R^n}\, |{\tilde \Lambda}_{21} (t,x)| \le \| f\|_{\gamma, d}
\, \sup_{x \in \R^n}  \,  \{ \E \big [ | E_0 (X_t^x - Y_t^x ) |^{
{ \gamma}} \, |J^2_{ij} (t,x) |\big ] \}  \,
\le c \| f\|_{\gamma, d} \; t^{\frac{ 2\gamma -2}{2 } \, - \, h \,
-h' }.
 \eeqr
 Using Lemma \ref{tec}, we infer  $\ds{ \sup_{x \in \R^n, \; t \in
(0,1)} \; | {\tilde \Lambda}_{22} (t,x)| = c_5 < \infty }$ and
 this gives the assertion.

%
%

\hh {\it - \  We prove  the   estimate  (iv) in \eqref{c7}, for $t
\in (0,1)$ and $\gamma \in (2,3)$ non-integer.}

\vskip 1mm \hh We have to show that, for any $h $, $0 \le h \le
k$,
 \beq \label{ty}
\sup_{x \in \R^n}\, \| P_t f (x +  \cdot ) \|_{C_b^{\gamma / 2h+1}
(E_h (\R^n))} \, \le \, c \| f\|_{\gamma, d}, \;\;\; f \in {\cal
C}^{\gamma}_d (\R^n),\;\; t \in (0,1).
 \eeq Fix  the integer $h$, $f \in
{\cal C}^{\gamma}_d (\R^n)$ and  consider $ \triangle_{v_h}^3 (P_t
f) (x) = P_tf(x) - 3 P_t f(x+ v_h) + 3 P_t f(x+ 2v_h) - $ $ P_t
f(x+ 3v_h)$, for $x \in \R^n$,  $v_h \in E_h (\R^n)$ with $|v_h|
\le 1$ and $v_h \not
 =0$. By \eqref{nag}
 the assertion \eqref{ty} is equivalent to the
 estimate
 \beq \label{ff}
 \sup_{x \in \R^n} \, |\triangle_{v_h}^3 (P_t f) (x)|  \le c_1 \, \|
 f\|_{\gamma, d} \, |v_h|^{\frac{\gamma}{2h+1}},\;\;\; t \in
 (0,1),
\eeq
 where $c_1$ is independent on $f$, $t$ and $v_h$. We prove
 \eqref{ff} considering first  the
  case of $|v_h| \le t^{\frac{2h+1}{2}}<1$ and then the
   case of   $1 \ge |v_h| >
  t^{\frac{2h+1}{2}}$  (compare with \cite[page 148]{Lu}).

\hh { (a)} \ Let $|v_h| \le t^{\frac{2h+1}{2}} <1$. Using the mean
value theorem  and  (iii) in \eqref{c7}, we get
 \beqr
   \sup_{x \in \R^n} \, |\triangle_{v_h}^3 (P_t f) (x)| \, \le
 \,  \sup_{x \in \R^n,\, i,\, j, \, r \in I_h} \,
  \| D_{ij r }^3
P_t f \|_0 \,  |v_h|^3
 \\
\le c \, \frac{1}{t^ { 3 h \,  + \, \frac{3 - \gamma}{2} } } \; \|
f\|_{\gamma, d} \, |v_h|^3  \, \le  \, c \, \frac{1}{|v_h|^ {
\frac{2}{2h+ 1} \,  \frac{6 h  + 3 - \gamma}{2} } } \; \|
f\|_{\gamma, d} \, |v_h|^3 \, =  \, c \, \|
 f\|_{\gamma, d} \, |v_h|^{\frac{\gamma}{2h+1}},\;\;\; t \in
 (0,1).
 \eeqr
{(b)} \ Let $1 \ge |v_h| > t^{\frac{2h+1}{2}}$. We first estimate
$ \nor
 e^{tA} v_h \nor $. To this purpose we use that
 \beq \label{ll}
 \| E_i e^{t  A} E_h
\|_L \le c  t ^{i - h }, \; \;\; 0 \le h \le i \le k; \;\;\; \|
E_i e^{t  A} E_h \|_L \le c  t, \;\; 0 \le i < h \le k, \;\; t \in
[0,1]
 \eeq
(see \cite[Lemma 3.1]{Lu})  where $c = c(A)>0. $
 Since
$ t \le |v_h|^{\frac{2}{2h + 1}} \le
 1$,  we get
   \beq \label{ci0}
    \bal
 \nor  e^{tA} v_h \nor = \sum_{i=0}^{h-1} |E_i e^{tA} E_h v_h
|^{\frac{1}{2 i +1}} + \sum_{i=h}^{k} |E_i e^{tA} E_h v_h
|^{\frac{1}{2 i +1}} \,
 \le c_1 \sum_{i=0}^{h-1} | t  v_h |^{\frac{1}{2 i +1}} \\
 + c_1
\sum_{i=h}^{k} t^{\frac{i-h} {2i +1}} | v_h |^{\frac{1}{2 i +1}}
\le c_1 \, h \, |v_h|^{1 / (2h+1)} + c_1 \sum_{i=h}^{k}  | v_h
|^{\frac{1}{2 i +1} +  \frac{i-h} {2i +1} \frac{2}{2h+1}}   \le
c_1 \, (k+1) \, |v_h|^{\frac{1}{2h+1}}.
 \eal \eeq
 To finish the proof we will use the Girsanov theorem, see \eqref{ls}.
 First  note that
 \beqr
\triangle_{v_h}^3 (P_t f) (x) \\  = \E \Big [ f(Z_t^x) \Phi (t,x)
 - 3f(Z_t^{x+ v_h}) \Phi (t,x+ v_h)
\\ + 3f(Z_t^{x+ 2v_h}) \Phi (t,x+ 2v_h) -
 f(Z_t^{x+ 3v_h}) \Phi (t,x+ 3v_h)
 \Big]
\\  =  A_1 (t,x) \, + \, A_2 (t,x), \;\;\; \mbox {where}
\\
 A_1(t,x) = \E \Big [ \Big ( f(Z_t^x)
 - 3f(Z_t^{\, x+ v_h})
+ 3f(Z_t^{\, x+ 2v_h})  -
 f(Z_t^{\, x+ 3v_h}) \Big ) \, \Phi (t,x)
 \Big],
\\
 A_2 (t,x) =  3 \E \big [ f(Z_t^{x+ v_h})  (\Phi (t,x) - \Phi
(t,x+ v_h) ) \big ]  + 3 \E \big [ f(Z_t^{x+ 2v_h})  (\Phi (t,x +
2 v_h) - \Phi (t,x ) ) \big ]
\\   +  \E \big [ f(Z_t^{x+ 3v_h})  (\Phi (t,x) - \Phi
(t,x+ 3v_h) ) \big ].
 \eeqr
Let us consider $A_1$. We find, for any $x \in \R^n$, $t \in
(0,1)$, thanks to Lemma \ref{zyg1},
 \beqr
&&  |A_1(t,x)|
\le \E \big [  \, |\triangle_{e^{tA}v_h}^3 f( e^{tA}x + Z_t^0)|
\Phi (t,x) \big ] \le  \| f\|_{\gamma, d} \, \nor e^{tA} v_h
{\nor}^{{\gamma}} \, \le  c  \| f\|_{\gamma, d} \, |
v_h|^{\frac{\gamma}{2h+1}}
 \eeqr
(in the last inequality we have used \eqref{ci0}). It remains to
treat $A_2$. We have:
 \beqr
 A_2 (t,x) &=& A_{21} (t,x) + A_{22} (t,x),
 \;\;\; \mbox {where}
 \\
 A_{21} (t,x) &=&  \E \Big [ f (Z_t^x) \Big ( \Phi (t, x) -  3
\Phi (t, x+ v_h) +  3 \Phi (t, x+ 2v_h)  -  \Phi (t, x+ 3 v_h)
\Big ) \Big ],
\\
 A_{22} (t,x) & =& 3 \E \big [ \big ( f(Z_t^{x+ v_h}) - f(Z_t^x)
\big )  \, \big ( \Phi (t,x) - \Phi (t,x+ v_h) \big ) \big ] \, \\
& + &\, 3 \E \big [ \big ( f(Z_t^{x+ 2v_h}) - f(Z_t^x) \big ) \,
\big( \Phi (t,x + 2 v_h) - \Phi (t,x ) \big ) \big ]
\\   &+ & \E \big [ \big ( f(Z_t^{x+ 3v_h}) - f(Z_t^x) \big ) \, \big(
 \Phi (t,x) - \Phi
(t,x+ 3v_h) \big) \big ].
 \eeqr
 In order to treat  $A_{21}$, remark that the map:
  $x  \mapsto \Phi (t,x) $ is three times
    Fr\'echet
differentiable from $\R^n$ with values in in $L^1 (\Omega)$.
 We need to estimate the norm of the first, second and third Fr\'echet
  derivatives of $\Phi (t,x)$; these Fr\'echet
  derivatives will be  indicated with
  $ D_x \Phi (t,x)$, $ D_{xx}^2 \Phi (t,x)$
  and $ D_{xxx}^3 \Phi (t,x)$ respectively.

   For any
$x,\, h \in \R^n$, we find (setting $G = Q^{-1/2}F$)
 \beqr
 D_x \Phi (t,x)[h] & = & \Phi (t,x) \int_0^t {\langle} DG(Z_s^x)
e^{sA}h , dL_s^x{\rangle}
 \\
&= &\Phi (t,x) \Big ( \int_0^t {\langle} DG(Z_s^x) e^{sA}h , dW_s
{\rangle} -   \int_0^t {\langle} DG(Z_s^x) e^{sA}h , G (Z_s^x )
{\rangle} ds \Big ),
  \eeqr
 since  $L^x_s : = W_s - \int_0^s G(Z^x_r)dr $, $ s \in
 [0,t]$
  (see \eqref{ls}).
  By  the Girsanov theorem, we have
 \beq \label{cosi}
\E  | D_x \Phi (t,x) [h] |
   = \E \Big | \int_0^t
{\langle} DG(X_s^x) e^{sA}h , dW_s {\rangle} \Big |
 \le \, e^{\| A\|_L} |h| \,  t^{1/2} \, \| DG \|_{0} \le  e^{\| A\|_L}
  |h| \, \| DG \|_{0}
 , \eeq
  for any $
  t \in [0,1], \; h \in \R^n.
 $ It follows that $\|D_x \Phi (t,x)  \|_{L (\R^n, L^{1}(\Omega))} \le
  e^{\| A\|} \| Q^{-1/2}_0\|_{L(\R^{{\tilde p}})} \| DF\|_0$, $t \in
  [0,1]$.
  Similarly, we have for the second Fr\'echet derivative
 \beqr
 D_{xx}^2 \Phi (t,x)[h][k] =
 \Phi (t,x)
 \Big (\int_0^t {\langle} DG(Z_s^x)
e^{sA}h , dL_s^x{\rangle} \Big ) \Big (\int_0^t {\langle}
DG(Z_s^x) e^{sA}k , dL_s^x{\rangle} \Big )
 \\
+ \, \Phi (t,x) \Big (  \int_0^t {\langle} D^2 G(Z_s^x) [e^{sA}k]
[e^{sA}h] , dL^x_s {\rangle} -   \int_0^t {\langle} DG(Z_s^x)
[e^{sA}h] , DG (Z_s^x )[e^{sA}k] {\rangle} ds \Big ),
  \eeqr
 $h, k
\in \R^n$. It follows, by the Girsanov theorem,
$$
\E  | D_{xx}^2 \Phi (t,x) [h] [k]| \,
    \le \, c_2 \,  |h| |k| \, (\| DG \|_{0}^2 + \| D^2 G \|_{0}),
      \;\; \mbox{ for any} \;
  t \in [0,1], \; h, k  \in \R^n.
$$
 In a similar way we get
 $$
\E  | D_{xxx}^3 \Phi (t,x) [h] [k] [u]| \,
    \le \, c |h| |k| |u| \, \big( \| DG \|_{0}^3 + \| DG \|_{0}^2 +
     \| D^2 G \|_{0}^2 +
     \| D^3 G \|_{0} \big) \le C_1  \, |h| |k| |u| ,
$$
 for any
  $t \in [0,1], \; h , k , u\in \R^n$,
  where $C_1 = C_1 (\| A\|_{L}\, , \nu_1, \tilde p ,
   \| DF \|_{0}, \| D^2F \|_{0},
  \| D^3 F \|_{0})>0$.
 Using the last estimate, we find
\beqr | A_{21} (t,x)| &\le& \| f\|_0  \, \| \Phi (t, x) - 3 \Phi
(t, x+ v_h) +  3 \Phi (t, x+ 2v_h)  -  \Phi (t, x+ 3 v_h) \|_{L^1
(\Omega)} \\ &\le & \| f\|_0 \sup_{|u|\le 1,  |h|\le 1, |k|\le 1,
\, x \in \R^n} \,
 \| D_{xxx}^3 \Phi (t,x) [h] [k] [u] \|_{L^1(\Omega)} \, \, |v_h|^3
\le \,  C_1  \| f\|_0 \, |v_h|^3,
 \eeqr
 $x\in \R^n,$ $ t \in [0,1].$ It remains to consider $A_{22}$. This
is the sum of three terms which can be treated in the same way.
Let us estimate the first term (without the factor $3$). By
\eqref{ci0}, we find
 (recall that $\gamma \in (2,3)$)
 \beqr
 \E \big | \big ( f(Z_t^{x+ v_h}) - f(Z_t^x) \big )  \, \big (
\Phi (t,x) - \Phi (t,x+ v_h) \big ) \big | \,
 \\
\le  \| f\|_{\gamma,d} \nor  e^{tA} v_h \nor \, \E  |\Phi (t,x) -
\Phi (t,x+ v_h)| \le  c  \| f\|_{\gamma,d}  \, |  v_h
|^{\frac{1}{2h +1}}
 \, \|\Phi (t,x) -
\Phi (t,x+ v_h) \|_{L^1 (\Omega)}.
 \eeqr
 By \eqref{cosi}, since  $|v_h| > t^{\frac{2h+1}{2}}$,
 \beqr
\E \big | \big ( f(Z_t^{x+ v_h}) - f(Z_t^x) \big )  \, \big ( \Phi
(t,x) - \Phi (t,x+ v_h) \big ) \big | \le
 e^{\|A \|_L} \| f\|_{\gamma,d}   |v_h|^{\frac{1}{2h +1}}  \, |v_h| \,  t^{1/2}
 \, \| DG \|_{0}
\\
 \le  c' \, |v_h|^{\frac{2}{2h+1} \, + 1} \,  \| f\|_{\gamma,d}.
 \eeqr
 We obtain
 $
  {\sup_{x \in \R^n} |A_{22} (t,x)| \le c_3
  \, |v_h|^{\frac{3 +  2h }{2h+1}} \,  \| f\|_{\gamma,d},\;\; }$ $t
  \in (0,1)$. Using the estimates for $A_1(t,x)$ and  $A_2(t,x)$,
    assertion  \eqref{ff} follows. This completes the proof.
\qed

\section {Elliptic and parabolic Schauder estimates
}

 Here we  prove elliptic and parabolic Schauder estimates for
${\cal A}$ using  the $L^{\infty}$-estimates of the previous
section. Our method is   different  with respect to
 \cite{Lu}, \cite{Ce}, \cite{Lo} and \cite{Sa} (see Theorems \ref{buon} and
  \ref{buon1}).
 Before proving Schauder estimates,
   we   show existence and
 uniqueness of distributional solutions for \eqref{1} and \eqref{2}.


\hh Let $\lambda >0$ and $f \in C_b(\R^n)$ (i.e., $f$ is uniformly
continuous and bounded on $\R^n$). We say that a function $u \in
C_b (\R^n)$ 
is a {\it distributional solution} to
the elliptic equation
 \beq \label{e3}
  \lambda u (x) - {\cal A} u(x) = f(x),\;\;\;\; x \in \R^n,
 \eeq
 if $ \lambda \int_{\R^n} u(x) \phi(x) dx = \int_{\R^n} u(x) {\cal A}^* \phi(x)
  dx + \int_{\R^n} f(x) \phi(x) dx $, for any $\phi \in C_0^{\infty}
  (\R^n)$, where ${\cal A}^*$ is the formal adjoint of ${\cal A}$,
  i.e.,
$$
{\cal A} ^*  \phi (x)= \frac{1}{2} \mbox{Tr} (Q D^2 \phi (x))  -
\lan Ax +
 F(x), D \phi (x) \ran - \phi (x)[ \mbox {div} F(x) + \mbox {Tr}(A)],
  \;\; x \in
 \R^n.
$$
\noindent Let $g \in C_b(\R^n)$, $T>0$ and $H : [0,T] \times \R^n
\to \R $ be a continuous and bounded function.
   We say that a continuous and bounded function $v : [0,T] \times
\R^n \to \R$ such that
 $v(0,x) = g(x)$, $x \in \R^n$, is a {\it
  space-distributional solution} to the parabolic Cauchy problem
 \beq \label{e4}
 \left\{ \bal \partial_t v(t,x) & = {\cal A} v
(t,x) \, +\, H(t,x), \;\;\; t \in (0,T],
 \; x \in \R^n,\\
      v(0, x) & = g(x),\;\; x \in \R^n.
\eal \right.
 \eeq
if the following conditions hold:

\hh (i) $v(t, \cdot ) \in C_b(\R^n)$ uniformly in $t \in [0,T]$;
(i.e.,
 for any $\epsilon >0$, there exists $\delta >0$ such that
 if $y \in \R^n$ and  $|y|< \delta$, we have
   $  \sup_{t \in [0,T],\, x \in \R^n} \;
    |v(t, x+y) - v(t,x)| <
 \epsilon )$.

\vskip 1mm \noindent (ii) for any test function $\phi \in
C^{\infty}_0 (\R^n)$,  the real mapping: $t \mapsto \int_{\R^n}
v(t,x) \phi (x) dx $ is continuously differentiable on $[0,T]$
and moreover
 \beq \label{fi}
  \bal
 \frac {d}{dt} \Big (\int_{\R^n} v(t,x) \phi (x) dx \Big )= \int_{\R^n}
  v(t,x) {\cal A}^* \phi (x) dx + \int_{\R^n}
  H(t,x)  \phi (x) dx, \;\; t \in [0,T].
  \eal \eeq
 \bth \label {eli}  Let $\lambda>0$ and $f \in C_b(\R^n)$. Then
 there exists a unique  distributional solution $u \in C_b(\R^n)$
 to the  equation
 \eqref{e3}. Moreover $u$ is given by
 \beq \label{lap}
 u(x) = \int_0^{\infty} e^{- \lambda t  } (P_t f) (x) dt
  = \int_0^{\infty} e^{- \lambda t  } P_t f (x) dt,\;\;\; x
 \in \R^n,
 \eeq
 where $P_t$ is the diffusion semigroup introduced in
 \eqref{pt}.

 \vskip 1mm
\noindent  Let $g \in C_b(\R^n)$, $T>0$ and $H : [0,T] \times \R^n
\to \R $ be continuous and bounded. Then there  exists a unique
 space-distributional solution $v$ to the  Cauchy problem \eqref{e4}.
Moreover, setting $\int_0^t P_{t-s} H(s,x) ds $ $:=  \int_0^t
P_{t-s} \big (H(s, \cdot ) \big ) (x) ds$, we have
 \beq \label{cau}
 v(t,x) =
 P_t g(x) + \int_0^t P_{t-s} H(s,x) ds , \;\;\; x
 \in \R^n, \; t \in [0,T].
 \eeq
 \eth
 \bpf
 $\underline{{ {\it Uniqueness}.}}$ \  We first consider the {\it elliptic
 case.}
 Fix  $\lambda >0 $ and let  $u \in C_b (\R^n)$ be any distributional
solution to \eqref{e3} with $f =0$.

 Take a function $\rho \in C^{\infty}_0 (\R^n)$ such that $\|
\rho\|_{L^1(\R^n) } =1$, $0 \le \rho \le 1 $ and   $\rho (x) =0$
 if $|x| \ge 1$. Define a sequence of mollifiers $(\rho_m )
  \subset C^{\infty}_0(\R^n)$, $\rho_m (x) : = m^n\rho(mx)$, $x \in
  \R^n$, $m \in \N$. Consider the functions $u_m \in C^{\infty}_b (\R^n)$
    obtained by
  convolution of $u$ with $\rho_m$, i.e., $u_m = u * \rho_m$.
 Setting $C(x) : = Ax + F(x)$, $x \in \R^n$, we use
the
   identity:  $${\cal A}^* [\rho_m (x - \cdot)](y) + \lan C(x) - C(y),
    D{\rho_m} (x-y)\ran  + \, \rho_m (x-y) \mbox{div}C(y) =
    {\cal A} [\rho_m ( \cdot - y)](x),
    $$ $ x, y \in \R^n,$ and  get
 \beq \label{fd}
 \bal
 {\cal A} u_m (x) = \int_{\R^n} u(y) \,  {\cal A} [\rho_m (\cdot \, - y
 )](x) \, dy   = \int_{\R^n} u(y) \,  {\cal A}^* [\rho_m (x  -
 \, \cdot
 )](y) \, dy   +  R_{m, 1} (x) + R_{m, 2} (x)
 \\
 = \int_{\R^n} \lambda u(y) \,   \rho_m (x \, - y
 ) \, dy  \, + \, R_{m, 1} (x) + R_{m, 2} (x),
 \; \;\;  \mbox{where}
\\  R_{m,1}(x) = \int_{\R^n} u(y) \, \mbox{div} C(y) \, \rho_m (x -y
) dy , \\   R_{m,2} (x) =  \int_{\R^n} u(y) \, \lan C (x) - C(y),
D \rho_m (x -y ) \ran  dy.
 \eal
 \eeq
 Changing variable as in   \cite[page 559]{Lo1} we obtain
 \beqr
 R_{m,2}(x) =
 m \int_{\R^n} u( x -  \frac{z}{m} ) \,
  \,  \lan C (x) - C ( x - \frac{z}{m}), D \rho (z) \ran
  \, dz.
 \eeqr
It follows that $R_{m,2} $ converges as $m \to \infty$, {\it
uniformly on $\R^n$,} to the function
$$ \ds{ x \, \mapsto \,
 u(x)   \sum_{i, k=1}^n
 \int_{\R^n} D_k C_i (x) z_k  D_i
\rho  (z) \, dz = - u(x) \mbox{div} C(x). }
$$
 On the other hand, it is easy to see that $R_{m,1}$ converges as $m \to \infty$,
 uniformly on $\R^n$, to $ u  \,  \mbox{div} C.$ It follows that
 $\lim _{m \to \infty} (R_{m,1} + R_{m,2}  ) =0$ in $C_b (\R^n)$.
 Hence we have  obtained
  $$
 \lim_{m \to \infty} \big ( \| {\cal A} u_m - \lambda u \| _0  \;
  + \; \|  u_m - u \| _0 )  \; =\; 0.
 $$
  By the classical maximum principle (see \cite{Lu3}) we deduce  that
 $ \ds{
 \| u_m \|_{0} \le \frac{1}{\lambda} \|  \lambda u_m - {\cal A} u_m
 \|_0}$. Letting  $ m   \to \infty$, we find that
  $\| u \|_0 =0$ and this gives the assertion.

\vskip 1mm \ We  prove now  uniqueness in the {\it parabolic
case}. To this purpose, we  take $H =0$ and $g=0$ in \eqref{e4}
 and consider any space-distributional solution $v$.
  We introduce as before a sequence of mollifiers $(\rho_m) \subset
  C^{\infty}_0 (\R^n)$ and define
$$
 v_m (t,x) = \int_{\R^n} v(t, y ) \rho_{m} (x-y) dy,\;\; t \in
 [0,T],\; x \in \R^n,\; m \in \N.
$$
It is clear that $v_m  $ is continuous and bounded on $[0,T]
\times \R^n$. Moreover, there exist continuous and bounded
  spatial  partial derivatives of  $v_m$ on $[0,T] \times \R^n
  $ of any  order. Thanks to
 assumption (i), $v_m$ converges to $v$ as $m \to \infty$
  {\it uniformly} on $[0,T] \times \R^n$.

We have,
 by \eqref{fi},
for $t \in [0,T]$, $x \in \R^n$,
 \beq \label{333}
 \bal
 {\partial }_t v_m (t,x) = \int_{\R^n} v(t,y) \,
  {\cal A}^* [\rho_m (x \, - \,
  \cdot
 )](y) \, dy
  = \int_{\R^n} v(t,y) \,  {\cal A} [\rho_m (\cdot  \, -
 \, y
 )](x) \, dy  \\  + \,\, S_{m, 1} (t,x) \, + \, S_{m, 2} (t,x)
  \;
 = \, {\cal A}v_m (t,x)  \, + \, S_{m, 1} (t,x) + S_{m, 2} (t,x),
 \; \;\;  \mbox{where} \;
 \\  S_{m,1}(t,x) = - \int_{\R^n} v(t,y) \, \mbox{div} C(y) \, \rho_m
(x -y ) dy , \\  S_{m,2} (t,x) =  -  \int_{\R^n} v(t,y) \, \lan C
(x) - C(y),   D \rho_m (x -y ) \ran \, dy.
 \eal
 \eeq
 Remark that $ { \lim_{m \to \infty}
  \, \sup_{t \in [0,T],\, x \in \R^n} | S_{m,1}(t,x) + S_{m,2}(t,x)|
 =0}$.
 Moreover, since $v_m$ is a classical solution to
$$
 \left\{ \bal \partial_t v_m(t,x) & = {\cal A} v_m
(t,x) \, +  S_{m,1}(t,x) + S_{m,2}(t,x),  \;\;\; t \in (0,T],
 \; x \in \R^n,\\
      v_m(0, x) & = 0,\;\; x \in \R^n.
\eal \right.
$$
 by the classical parabolic maximum principle (see
 \cite[Chapter 8]{Kry}) we have
 $$ \ds{ \sup_{t \in [0,T],\, x \in \R^n} | v_m (t,x)|
  \; \le \;  T
   \sup_{t \in [0,T],\, x \in \R^n} | S_{m,1}(t,x) + S_{m,2}(t,x)|
 }.$$
 Letting $m \to \infty$ we obtain that $v =0$ and this proves
  the assertion.

\hh $\underline{{ {\it Existence}.}}$ \  We first consider the
 {\it elliptic case}  and   prove that $u$ given in \eqref{lap}
 is the distributional solution.
  It is clear that    $u \in C_b (\R^n)$. In the following
computations we will use that there exists the classical partial
 derivative $\partial_t (P_t f)(x)$, for $t>0$ and $x \in \R^n $,
 and  $\partial_t (P_t f)(x)=
 {\cal A}(P_t  f )(x)$,  see \cite[Section
 4]{Pr}.

 By Corollary \ref{c3} we deduce that, for any
$M>0$, there exists $C_M >0$ such that
 \beq \label{vi} \sup_{|x| \le M}
 | {\cal A} (P_t f)(x) | \, \le \, {C_M}\, ({t^{- (1+k)}} +  1) \, \| f\|_0,
 \;\;\;
 t>0, \; f \in C_b (\R^n).
\eeq  We obtain, for any
 $\phi \in C^{\infty}_0 (\R^n)$, applying
 the Fubini theorem,
 \beqr
 \int_{\R^n}  u(x) \, {\cal A}^* \phi (x) =
 \int_0^{\infty} e^{- \lambda t  } dt \int_{\R^n} {\cal A} P_t f (x) \phi(x)
 dx = \lim_{\epsilon \to 0^+}
 \int_{\epsilon}^{\infty} e^{- \lambda t  } dt
  \int_{\R^n} {\cal A} P_t f (x) \phi (x)
 dx
\\
 = \lim_{\epsilon \to 0^+}
 \int_{\epsilon}^{\infty} e^{- \lambda t  } dt \int_{\R^n}
 \partial_t
P_t f (x) \phi(x)
 dx =
 \lim_{\epsilon \to 0^+} \int_{\R^n} \phi(x)dx \,
\int_{\epsilon}^{\infty} e^{- \lambda t }
\partial_t P_t
f (x)  dt \\
 = \lim_{\epsilon \to 0^+}
  \int_{\R^n} \Big ( - e^{- \lambda \epsilon} P_{\epsilon} f(x)  \, + \,
 \lambda \int_{\epsilon}^{\infty} e^{- \lambda t} P_t f(x ) dt \, \Big )
  \phi(x) = \int_{\R^n}  (- f(x) + \lambda u(x) ) \, \phi(x)
 dx.
 \eeqr
 We  deal now with  the
 {\it parabolic  case}  and   show  that $v$ given in \eqref{cau}
 is the space-distributional solution.
  We write
  \beq  \label{vf} v = v_1 + v_2, \;\;\; \mbox{
  where}
  \;\;  v_1 (t,x) = P_t g (x), \;\;\; v_2 (t,x) = \int_0^t P_{t-s} H(s,x)
  ds,
  \eeq
  $v_2 (0, \cdot ) =0$ ($v_1$ and $v_2$ are associated  to \eqref{cau} when $H=0$ and $g=0$
   respectively). First we deal with  $v_1$.
In \cite[Section 4]{Pr} it is verified that $v_1$ is  a continuous
and bounded function on $[0, \infty) \times \R^n$.  Moreover,
denoting by $\omega_g$ the modulus of continuity of $g$, we have,
for any $t \in [0,T]$, $x, y \in \R^n,$
 $| P_t g(x) - P_t g(y) | \le$
 $\E [\omega_g (|X_t^x -X_t^y|)] \le \omega_g (|x - y|\, e^{TL}),$
 where $L = \| A\|_{L}\,  +  \| DF\|_0$. This shows that $ v_1 (t,
\cdot ) \in C_b(\R^n)$, uniformly in $t \in [0,T]$.

\vskip 1mm  \noindent Since it holds  (in a classical sense) $\ds
{
\partial_t (P_t f) (x) = {\cal A}(P_t f) (x)}$, $t>0$, $x \in
\R^n$,   we have that  $t \mapsto \int_{\R^n} v_1 (t, x) \phi(x)
dx $ belongs to $C^{1}([0,T])$ and verifies \eqref{fi} (with
$H=0$).

\vv

\noindent  Let us treat $v_2$. By the first estimate in
\eqref{chi} we
 deduce, for any $f : \R^n \to \R $ continuous and bounded, for
 any $h \in \{0, \ldots, k \}$,
$$
 \| P_t f(x+ \cdot) \|_{C^{\frac{1}{2k+1}}_b (E_h (\R^n))} \le
 \| f\|_0^{1 - \frac{1}{2k+1}} \,\,  \|
  P_t f(x+ \cdot) \|_{C^{1}_b (E_h (\R^n))}^{\frac{1}{2k +1}}
 \le C \, t^{- 1/2} \| f\|_0, \; t \in (0,T],
$$
 $x \in \R^n,$ where $C$ is independent on $t$, $x$ and $f$. It follows that,
 for any $x,y \in \R^n$, $t \in [0,T]$,
 $$
|v_2(t,x) - v_2(t,y)| \le \int_0^t \frac{C}{ (t-s)^{1/2}}
 ds \, \sum_{h=0}^k |E_h (x-y)|^{\frac{1}{2k +1}} \,  \le
 c' \sqrt{T} |x-y|^{\frac{1}{2k +1}}.
 $$
 This shows that $v_2(t, \cdot ) \in C_b (\R^n)$, uniformly in
 $t$. Thanks to this  property, in order to verify that
 $v_2$ is continuous on $[0,T] \times \R^n$, it is enough to check
  that for any fixed $x \in \R^n$, $v_2(\cdot, x)$ is continuous on
  $[0,T]$.
 Since the continuity of $v_2 ( \cdot ,x)$ in $t=0$ is clear, we
consider continuity at  $t \in (0, T]$. We  write, for   $h $
sufficiently small,
  \beq \label{fu}
 v_2( t+h,x)- v_2(t,x) \, =\,
 \int_0^T \big [P_{ t +h -s}\, H(s,x) \, -\, P_{
t-s}H (s,x) \big ] ds.
 \eeq
(we have extended $P_t$ to negative values, setting
$P_{\eta}=0,\;\; \eta <0 $). By the dominated convergence theorem
 one deduces that
  $\lim_{h \to 0} v_2(t+h, x) = v_2(t,x)$. Thus $v_2$ is
  continuous on $[0,T ] \times \R^n$ and $v_2 (0, \cdot) = 0$. The
   boundedness of $v_2$ is clear.

\vv It remains to verify that $v_2$ satisfy \eqref{fi}. To this
purpose, we fix
 $t \in (0,T]$, $x \in \R^n$,
 and consider for $h>0$,
   see
 also \cite[pages 58-59 ]{Pr0},
$$
 \bal
 \frac{v_2(t+h, x) - v_2(t,x)}{h} = \Gamma_1(t, h,x) + \Gamma_2
(t, h, x),
\\ \Gamma_1(t, h,x) =   \frac{1}{h}
\int_t^{ t+h} P_{ t +h -s} H (s,x) ds, \;\;   \Gamma_2(t, h,x) \,
=\,
  \int_0^t \Big ( \frac { P_{ t+h -s}
-\,  P_{ t-s} }{h} \Big)  H(s,x)  ds. \eal $$
 We have: $
|\Gamma_1 (t, h,x) - H(t,x)|  \le
  \int_0^{ 1}  \E |
H (t + h - sh , X_{sh}^x) \,  - \, H(t,x) | ds
\rightarrow \; 0 \;  \mbox{as}$ $h$ tends to $0^+$,
 by the dominated convergence theorem. It follows that,
   for any  $\phi \in C^{\infty}_0(\R^n) $,  $\lim_{h \to
0^+} $ $ \quad $ $\; \int_{\R^n} \Gamma_1(t, h,x) \phi(x) dx$ $ =
\int_{\R^n} H(t,x) \phi(x) dx  $.

Concerning $\Gamma_2$, we first note that, thanks to \eqref{vi},
for any $t >s \ge 0$,
 \beq \bal & \lim_{h \to 0^+} \int_{\R^n} \! \! \! \big ( \frac { P_{ t+h
-s} H(s,x) -  P_{ t-s} H(s,x)}{h}  \big ) \phi(x) dx  =
\int_{\R^n} {\cal A} [ P_{ t-s} H(s, \cdot)](x) \,
  \phi(x) dx .
\\
& \mbox{ By the Fubini theorem we get} \;\;
 \lim_{h \to 0^+}
\int_{\R^n} \Gamma_2 (t,h,x) \phi(x) =
\\ & = \int_0^t ds \int_{\R^n}  P_{ t-s} H(s, x)   {\cal A}^* \phi(x)
dx = \int_{\R^n }  {\cal A^*}
 \phi (x) dx
 \int_0^t
P_{ t-s} H(s, x)  ds , \; t \in ]0,T].
 \eal \eeq
 It follows easily that  the map $t \mapsto \int_{\R^n} v_2(t,x)
\phi (x)$ belongs to $C^1 ([0,T])$ and   verifies \eqref{fi} (with
$g=0$) for $t \in [0,T]$.  This finishes the proof.
 \epf
 The next theorems provide elliptic and parabolic  Schauder estimates.
\bth \label{buon}
 Let $\theta \in (0,1)$ and  $\lambda >0$. For  any
 $f \in {\cal C}^{\theta}_d (\R^n)$
 there exists a unique distributional solution  to the elliptic equation
  \eqref{e3}.  Moreover $u \in C^{2 + \theta}_d (\R^n)$ and there exists
  $c = c ( \lambda, \theta, \nu_1,  \nu_2, A, \tilde p,
    n, \| D^{} F\|_0, \| D^{2} F\|_0
  \| D^{3} F\|_0 )$,  such that
 \beq \label{e2}
 \| u \|_{{2 + \theta}, d} \, \le \,
  c \| f\|_{\theta, d}.
 \eeq
 \eth
 \bpf  Uniqueness follows by Theorem \ref{eli}.
   We need to investigate the  regularity properties of the
   function
   $u \in C_b (\R^n)$ given in \eqref{lap}.

   We first prove that $u  (z +  \cdot ) \in C^{2 + \theta}_b
 (E_0 (\R^n))
  $, for any $z \in \R^n$,  and
 \beq \label{dd}
 \sup_{z \in \R^n} \| u  (z +  \cdot ) \|_{C^{2 + \theta}_b
 (E_0 (\R^n))} \, \le \, C \, \| f\|_{\theta, d}.
 \eeq
 It is clear by  the estimates
  \eqref{elli1} that there exist the partial derivatives
   $D_i u $ and $D_{ij}^2 u $ on $\R^n$, for any $i,j \in I_0$.
    Moreover $D_i u $ and $D_{ij}^2 u $  are continuous
   and bounded  on $\R^n$ and
   $\| D_i u\|_0 $ $ + \,
    \| D_{ij}^2 u\|_0 \le$ $ c \| f\|_{\theta, d}$.

We will prove now that  $D_{ij}^2 u \in {\cal C}^{\theta}_d
(\R^n)$ when $i,j \in I_0$. This will imply \eqref{dd}. To this
purpose,  we fix $v_h \in E_h (\R^n)$, for  $0 \le h \le k$, with
 $|v_h|\le 1$, and
compute
 \beq \label{vai}
 \bal | D_{ij}^2 u(x + v_h) - D_{ij}^2 u (x)|  &\le
 \int_0^{\infty} e^{-
\lambda t } | D_{ij}^2 P_t f (x+ v_h) -  D_{ij}^2 P_t f (x)| \, dt
\; =
 \; u_1(x) +  u_2(x), \\
u_1(x) &= \int_0^{|v_h|^{\frac{2}{2h+1}}} e^{- \lambda t } |
D_{ij}^2 P_t f (x+ v_h) - D_{ij}^2 P_t f (x)| dt; \\  u_2(x) &=
\int_{|v_h|^{\frac{2}{2h+1}}}^{\infty} e^{- \lambda t } | D_{ij}^2
P_t f (x+ v_h) - D_{ij}^2 P_t f (x)| \, dt,\;\; x \in \R^n. \eal
 \eeq
 In order to estimate $u_1(x)$ we use  (b)
 in \eqref{elli1}. We find
$$
u_1(x) \le c  \| f\|_{\theta, d} \,
\int_0^{|v_h|^{\frac{2}{2h+1}}} t^{\frac{\theta}{2} - 1} \, dt \le
C \| f\|_{\theta, d} \, \, |v_h|^{\frac{ \theta}{2h+1}}.
$$
 Concerning $u_2(x)$ we use estimate (c) in \eqref{elli1}. This
 gives
$$ |D_{ij}^2 P_t f (x+ v_h) - D_{ij}^2 P_t f (x)| \le
 |v_h| \sup_{r \in I_h}
  \| D_{ij r }^3 P_t f \|_0 \le c \| f\|_{\theta, d} \, \Big (\frac{1}{t^ {
\frac{3 - \theta}{2} + h }  } + 1 \Big )\, |v_h|, \; \; t>0.
$$
We get
 \beqr
 u_2(x) \le  c  \| f\|_{\theta, d} \, |v_h|\,
\int_{|v_h|^{\frac{2}{2h+1}}} ^{\infty} e^{-\lambda t} \big (
 t^{\frac{\theta}{2} - \frac{3}{2} - h } \, + 1 \,  \big ) dt \,
  \le c' \, \big(\frac{|v_h|}{\lambda} +  |v_h|^{\frac{ \theta}{2h+1}} \big)
  \| f\|_{\theta, d}
\\  \le C_1 \| f\|_{\theta, d} \,  |v_h|^{\frac{ \theta}{2h+1}},
 \;\;\; x \in \R^n. \eeqr
 It follows that $| D_{ij}^2 u(x + v_h) - D_{ij}^2 u (x)|  \le C \| f\|_{\theta, d}
  \, \, |v_h|^{\frac{ \theta}{2h+1}}$ and so
   \eqref{dd} is  proved.

   We verify  that $u  (z +  \cdot ) \in C^{\frac{ 2 + \theta}{2h
 +1}}_b
 (E_h (\R^n))
  $, for any $1 \le h \le k$, and moreover
 \beq \label{dd2}
 \sup_{z \in \R^n} \| u  (z +  \cdot ) \|_{C^{\frac{2 + \theta}{2h +
  1}}_b
 (E_h (\R^n))} \, \le \, C \, \| f\|_{\theta, d}.
 \eeq
We fix   $v_h \in E_h (\R^n)$, for  $1 \le h \le k$,
 with $|v_h|\le 1$, and
compute
 \beq \label{rit1}
 \bal |  u(x + v_h) -  u (x)|  \le
 \int_0^{\infty} e^{-
\lambda t } |  P_t f (x+ v_h) -   P_t f (x)| \, dt \; =
 \; u_1(x) +  u_2(x), \;\;\; \mbox{where}
\\
u_1(x) = \int_0^{|v_h|^{\frac{2}{2h+1}}} e^{- \lambda t } |  P_t f
(x+ v_h) - P_t f (x)| dt;  \\    u_2(x) =
\int_{|v_h|^{\frac{2}{2h+1}}}^{\infty} e^{- \lambda t } |  P_t f
(x+ v_h) -  P_t f (x)| \, dt,\;\; x \in \R^n. \eal
 \eeq
 In order to estimate $u_1(x)$ we use  (d)
 in \eqref{elli1}. We find
$$
u_1(x) \le c  \| f\|_{\theta, d} \,
 |v_h|^{\frac{ \theta }{2h+1}}
  \int_0^{|v_h|^{\frac{2}{2h+1}}}  dt \le C \| f\|_{\theta, d} \, \,
|v_h|^{\frac{ 2 + \theta}{2h+1}}.
$$
 Concerning $u_2(x)$ we use estimate (a) in \eqref{elli1}. We get
  (recall that  $h \ge 1$)
 \beqr
u_2(x) \le  c  \| f\|_{\theta, d} \, |v_h|\,
\int_{|v_h|^{\frac{2}{2h+1}}} ^{\infty} e^{-\lambda t} \big (
 t^{\frac{\theta}{2} - \frac{1}{2} - h } \, + 1 \,  \big ) dt \,
 \, \le \,
 C_1 \| f\|_{\theta, d} \, \, |v_h|^{\frac{ 2+ \theta}{2h+1}}
\eeqr and \eqref{dd2} follows. The proof is complete.
 \epf
\bth \label{buon1}
 Let $\theta \in (0,1)$, $T>0$, $g \in C^{2 + \theta}_d
  (\R^n)$  and let $H : [0,T] \times \R^n \to \R $
be a continuous  function such that ${ \sup_{t \in [0,T]}
 \| H(t, \cdot)\|_{\theta, d} < \infty}$.

  Then the Cauchy problem \eqref{e4} has a unique space-distributional
  solution $v$ such that  $v (t, \cdot ) \in C^{2 + \theta}_d
  (\R^n)$, $t \in [0,T]$.
   Moreover,  $D_i v$  and $D^2_{ij} v$ are continuous
  on $
   [0, T] \times \R^n$, for  $i,j \in I_0$, and
    there exists
  $c = c ( T, \theta, \nu_1, \nu_2, A, \tilde p,  n, \| D^{} F\|_0,$ $
   \| D^{2} F\|_0
  \| D^{3} F\|_0 )$,  such that
  \beq \label{e23}
 \sup_{t \in [0,T]} \| v(t, \cdot) \|_{{2 + \theta}, d} \, \le \,
  c \big ( \| g\|_{2+ \theta, d}\, +\,  \sup_{t \in [0,T]}
 \| H(t, \cdot)\|_{\theta, d} \big).
 \eeq
 \eth
 \bpf  Uniqueness follows by Theorem \ref{eli}.
 To prove the result,   we need to investigate the  space-regularity
  of the
   function
   $v  $ given in \eqref{cau}; we write $v = v_1 + v_2$ as in   \eqref{vf}.

  Concerning the function $v_1 =  P_t g$ the  estimate (iv) in \eqref{c7}
  with $\gamma = 2 + \theta$ gives immediately  \eqref{e23} with
  $v$ replaced by $v_1$ and $H=0$.
 In order to treat $v_2$,
$$
 v_2(t,x) = \int_0^t \E [H (s, X_{t-s}^x)]ds
 = \int_0^t \E [H (t - s, X_{s}^x)]ds,\;\; t \in [0,T],\, \,
  x \in \R^n,
$$
  we proceed as in the proof of Theorem
 \ref{buon}. 
  To this purpose,   set
  $ {\| H \|_{T, \theta} =
  \sup_{t \in [0,T]}\|H (t, \cdot )\|_{\theta, d}}$.
  We first prove that $v_2  ( t, z +  \cdot ) \in C^{2 + \theta}_b
 (E_0 (\R^n))
  $, for $t \in [0,T]$ and $z \in \R^n$,  and that
 \beq \label{dd1}
 \sup_{t \in [0,T],\, z \in \R^n} \| v_2  (t, z +  \cdot ) \|_{C^{2 + \theta}_b
 (E_0 (\R^n))} \, \le \, C \, \| H\|_{T, \theta}.
 \eeq
  It is clear by  the estimates
  \eqref{elli1} that there exist the spatial partial derivatives
   $D_i v_2 $   and $D_{ij}^2  v_2 $ on $ [0,T] \times \R^n $,
    for any $i,j \in I_0$.
    Moreover $D_i v_2 (t, \cdot) $ and $D_{ij}^2 v_2 (t, \cdot) $
    are  continuous   and bounded
    on $\R^n$ and  $\| D_i v_2 (t, \cdot)\|_0 $ $ + \,
    \| D_{ij}^2 v_2 (t, \cdot)\|_0 \le $ $  c \, \| H\|_{T, \theta}$, for any
     $t \in [0,T]$.

\vv To prove  assertion \eqref{dd1},   we fix $v_h \in E_h
(\R^n)$, for  $0 \le h \le k$,
 with $|v_h|\le 1$, and
compute as in \eqref{vai}
 \beqr
  | D_{ij}^2 v_2 (t, x + v_h) - D_{ij}^2 v_2 (t, x)|  \le
 \int_0^{t}  | D_{ij}^2 P_s H (t-s, x+ v_h) -
  D_{ij}^2 P_s H(t- s, x)| \, ds
 \\
\le  c \| H\|_{T, \theta} \, \int_0^{t \, \wedge  \,
|v_h|^{\frac{2}{2h+1}} }  s^{\frac{\theta}{2} - 1} \, ds \; + \;
c_1 \| H\|_{T, \theta} \, |v_h|\, \int_{t \, \wedge  \,
|v_h|^{\frac{2}{2h+1}} }^t \,   s^{\frac{\theta}{2} - \frac{3}{2}
- h }   \, ds
   \le
 \, c' \,  \| H\|_{T, \theta} \,  |v_h|^{\frac{ \theta}{2h+1}}
 \eeqr
($a \wedge b =$  min$(a,b)$)  and so the assertion
   \eqref{dd1} is  proved. In order to
    verify  that $v_2  (t, z +  \cdot ) \in C^{\frac{ 2 + \theta}{2h
 +1}}_b
 (E_h (\R^n))
  $, for any $1 \le h \le k$,
 $t \in [0,T]$,
   and that
 $
  \ds {\sup_{z \in \R^n,\, t \in [0,T]} \,
  \| v_2  (t, z +  \cdot ) \|_{C^{\frac{2 + \theta}{2h +
  1}}_b
 (E_h (\R^n))} } $ $\, \le \, C \, \| H\|_{T, \theta},
 $
 we proceed as in \eqref{rit1}.

In order   to prove
    the continuity
  of $D_i v$ and $D^2_{ij}v $ on $[0, T] \times \R^n$,
  $i, j \in I_0$,
 it is enough to show that,  for any fixed $x \in \R^n$,
 $D_i v(\cdot, x)$  and $D_{ij}^2 v (\cdot, x)$ are continuous on
  $[0,T]$. To this purpose, we write $x = x_0 + x_1$, where $x_0 =E_0 x$
   and $x_1 = x - E_0 x$, and  consider the   closed euclidean  ball $K$
   centered in $x_0$
   with radius 1. We already now that $ \| v(t, x_1 + \cdot )\|_{C^{2 + \theta}
    (K)} \le C_T$, for any $t \in [0,T]$. Using the continuity
     of $v$ on $[0,T] \times \R^n$ and a standard compactness
    argument we obtain  the assertion. Note that in particular
 $ \lim_{t \to 0^+} D_i v(t, x) = D_i g(x)$  and
  $\lim _{t \to0^+} D_{ij}^2 v (t, x) = D_{ij}^2 g(x)$, $x \in \R^n$.
 \epf

\section {Schauder estimates with variables coefficients $(q_{ij})$}

Here we consider a generalization of the operator $\cal A$, namely
 we deal with
  the operator $\tilde A$ in which  the diffusion matrix $Q$ depends
continuously on $x$, i.e.,
 \beq \label{opp}
 \tilde {\cal A} u (x) =
  \frac{1}{2} {\mbox{\rm Tr }}  (Q (x) D^2 u (x)) +
  {\langle}  Ax, D u (x) {\rangle}
 + {\langle}  F(x), D u (x){\rangle} , \;\;\;
 \; x \in \R^n.
 \eeq
Using a standard approach based on maximum principle, a priori
estimates and continuity method  (compare with \cite[Section
6]{Lu}) we will extend elliptic and parabolic Schauder estimates
of Section 4 to the operator $\aa$.

\begin{hypothesis} \label{hyy} {\em (i) \ there exists $\nu>0$ and an integer
 ${\tilde p}$, $1 \le \tilde p \le n$,  such that the symmetric matrix
  $Q (x)= (q_{ij} (x))_{i,j=1, \ldots ,n}$
   has the form
 \beq \label{q00} Q(x) = \left (
\begin{matrix} Q_0 (x) & 0 \\ 0 & 0
  \end{matrix} \right ), \;  \;\;\; x \in \R^n,
 \eeq
where $Q_0 (x)$  is a
    positive definite  ${\tilde p} \times {\tilde p}\,$-matrix such that
\begin{equation}
 \mbox{$  \nu \sum_{i =1}^{\tilde p}  \xi_i^2  \le
  \sum_{i,j =1}^{\tilde p} q_{ij} (x) \xi_i \xi_j \le \frac{1}{\nu}
  \sum_{i =1}^{\tilde p}
  \xi_i^2, \;\;\; \xi = (\xi_i)\in \R^p,\; x \in \R^n. $}
   \label{q01}
\end{equation}
 (ii) \ the vector field $F: \R^n \to \R^n$ satisfies (ii) and
 (iii) in  Hypothesis \ref{hy}.

\hh (iii) \ assumption (iv) in  Hypothesis \ref{hy} holds.

\hh (iv) \ There exists $\theta \in (0,1)$ such that $q_{ij} \in
{\cal C}^{\theta}_d (\R^n)$, for $i, j \in \{ 1, \ldots , {\tilde
p} \}$, and moreover there exists the
 limit
\begin{equation} \label{q02}
\lim_{|x|\mapsto \infty }Q_0(x) = Q_{0}^{\infty } \;\;  \mbox{in}
\; \, L(\Bbb R^{{\tilde p}}).
\end{equation}
}
\end{hypothesis}
\noindent Let us comment on these assumptions. Note that, for
every $x_0\in \Bbb R^n$, the operator with frozen second order
coefficients
\begin{equation}{\cal A} (x_0) =
\frac{1}{2}\mbox{\rm Tr} \,
 (Q(x_0)D^2\cdot
) + \langle F(x) + A x,D\cdot \rangle \label{841}
\end{equation}
 verifies Hypothesis \ref{hy} and therefore Theorems \ref{buon}
 and \ref{buon1} holds for  ${\cal A}( x_0)$. The same happens for
 the operator ${\cal A}^{\infty}  $ defined as in \eqref{841} but
 with $Q(x_0)$ replaced by  $Q^{\infty}$
  ($Q^{\infty }$ is the $n\times n$ matrix having $Q_{0}^{\infty }$
in the first $ {\tilde p} \times {\tilde p}$ block, and zero
entries in the other blocks;
 clearly    its  coefficients $q_{ij}^{\infty}$ verify
 \eqref{q01}).

To prove the next theorems it is crucial to remark
that the constants
  in the elliptic and
  parabolic Schauder estimates involving  ${\cal A} (x_0)$
  {\it does not depend on $x_0 \in \R^n$.}
 \begin{theorem} \label{t1}
 Consider the operator $\aa$  in \eqref{opp} under   Hypothesis \ref{hyy}.
  Then, for every $\lambda
>0$ and $f\in {\cal C}^{\theta }_{d}(\Bbb R^n)$ the elliptic problem
\begin{equation}
    \lambda u -  \tilde{{\cal A}}u =f
\label{ell8}
\end{equation}
has  a unique  solution $u\in {\cal C}^{2+\theta }_{d}(\Bbb R^n)$
 (here the first order term $\langle  A x , Du (x)\rangle $ is understood in
 distributional sense).
 Moreover there is $c>0$, independent of $f$ and $u$, such that
 Schauder estimates \eqref{e2} hold for \eqref{ell8}.
\end{theorem}
 \bpf
 We will only  sketch the proof which is not difficult.
 One  needs first a  {\it maximum principle} for \eqref{ell8}. We
 explain how this result can be obtained  arguing as in the proof of Theorem
 \ref{eli}. We write $\aa = {\cal A}_1 $ $+ \; {\cal A}_2$,
 where
  \beq
 {\cal A}_1 = \frac{1}{2}\mbox{\rm Tr} \,
 (Q(x)D^2\cdot
) \;\;\; \mbox{and} \;\;\;\;  {\cal A}_2 =  \langle F(x) + A x, D
\cdot \rangle. \label{84} \eeq
 Take any $u \in {\cal C}^{2 + \theta}_d (\R^n)$ which solves
\eqref{ell8}. Consider a sequence of mollifiers $(\rho_m)$ and set
 $u_m = u * \rho_{m}$;
we get, similarly to \eqref{fd},
$$
 \aa u_m (x) =  \int_{\R^n} {\cal A}_1 u( x - y)
 \rho_m (y ) dy + \int_{\R^n} u(y) \,  {\cal A}_2^* [\rho_m (x  -
 \, \cdot
 )](y) \, dy   +  R_{m, 1} (x) + R_{m, 2} (x),
$$
 $x \in \R^n$, $m \in \N$,
  where ${\cal A}_2^*$ is the formal adjoint of ${\cal A}_2$.
 One finds that $\aa u_m $ converges in $C_b (\R^n)$ to
  $\aa u$ as $m \to \infty$.
     By the classical maximum principle (see \cite{Lu3}) we deduce  that
 $
 \| u_m \|_{0} \le \frac{1}{\lambda} \|  \lambda u_m - \aa u_m
 \|_0$. Letting  $ m   \to \infty$, we find
 $\| u \|_{0} \le \frac{1}{\lambda} \|  \lambda u  - \aa u
 \|_0.$

 \hh \noindent {\it A priori estimates}  for \eqref{ell8} can be proved exactly as in
 the proof of \cite[Theorem 8.1]{Lu}. One assumes that $u \in {\cal C}
 ^{2 + \theta}_d (\R^n)$ is a solution to \eqref{ell8} and then by
 using a
  localization argument and the maximum principle
   one finds that there exists $C = C>0$ (independent
  on $f$ and $u$) such that
$$
 \| u\|_{2 +  \theta, d} \le C \| f\|_{\theta, d}.
$$
The {\it continuity method} allows to conclude the proof. For any
 $\epsilon \in [0, 1]$ one considers the problem
   \beq \label{cie1}
 \lambda u - (1- \epsilon ){\cal A}^{\infty} u - \epsilon
  \aa u = f,
 \eeq
 where $ (1- \epsilon ){\cal A}^{\infty} u (x)  +  \epsilon
  \aa u (x)= \frac{1}{2}\mbox{\rm Tr} \,
  \big ( [(1- \epsilon) Q^{\infty} + \epsilon  Q(x)] D^2 u (x)
\big ) +  $ $\langle F(x) + A x,D u (x) \rangle  $.

 Using the previous a priori estimates, it is straightforward
  to verify that
 the set of all $\epsilon$'s such that \eqref{cie1} is uniquely
 solvable in ${\cal C}^{2+  \theta}_d (\R^n)$ is non-empty, closed
 and open in $[0,1]$.
Taking $\epsilon =1$ in  \eqref{cie1} one
 finishes the proof.
\epf \noindent In order to state and prove Schauder estimates for
the parabolic Cauchy  problem involving $\aa$,  we define the
space ${\cal C}^{\gamma}_{T, d} $, $\gamma \in (0,3)$ non-integer.
This consists of all continuous functions $v : [0,T ] \times \R^n
\to \R$ such that
 $v(t, \cdot ) \in {\cal C}^{\gamma}_d (\R^n)$, $t \in [0,T]$,
  and moreover $\sup_{t \in [0,T] } \| v (t, \cdot)\|_{{\cal C}^{\gamma}_d
  (\R^n)} < + \infty$.
  $ {\cal C}^{\gamma}_{T, d}$ is a Banach space endowed with the
 norm  $\| \cdot  \|_{\gamma, T, d}$,
$$
 \| v \|_{\gamma, T, d} = \sup_{t \in [0,T] } \| v (t, \cdot)\|_{{\cal C}^{\gamma}_d
  (\R^n)}, \;\;\;\; v \in {\cal C}^{\gamma}_{T, d}.
$$
A function $v \in {\cal C}^{2+  \theta}_{T, d} $, $\theta \in
(0,1)$, solves the Cauchy  problem \eqref{e4} for  $\aa$ if $v(0,
x) = g(x)$, $x \in
  \R^n, $ and, for any $\phi \in C^{\infty}_0 (\R^n)$,
   the real mapping: $t \mapsto \int_{\R^n} v(t,x) \phi (x)
 dx $ is continuously differentiable on $[0,T]$ and
 verifies, for any $t \in [0,T]$ (see
 \eqref{84}$)$,
  \beq \label{cie}
\frac {d}{dt} \Big (\int_{\R^n} v(t,x) \phi (x) dx \Big )=
\int_{\R^n}   {\cal A}_1 v(t,  x) \, \phi (x) dx + \int_{\R^n}
v(t,  x) \, {\cal A}_2^* \phi (x) dx
 + \int_{\R^n}
  H(t,x)  \phi (x) dx.
\eeq
 \begin{theorem}
 Consider the operator $\aa$  in \eqref{opp} under   Hypothesis \ref{hyy}.
 Let $T>0$, $g \in {\cal C}^{2 + \theta}_d
 (\R^n)$ and $H \in {\cal C}^{\theta}_{T, d}$.
  Then  there exists a unique solution $v \in {\cal C}^{2 +
   \theta}_{T, d} $
  to the Cauchy  problem \eqref{e4} for  $\aa$.
   Moreover the spatial partial derivatives
    $D_i v$ and $D_{ij}^2 v $ are continuous on
    $[0,T] \times \R^n$, for $i, j \in I_0$, and
 there exists $c>0$, independent of $g$, $H$ and $v$, such that
 \beq \label{df}
  \|v \|_{2 + \theta, T, d}  \, \le \,
  c  \big (\| g\|_{2 + \theta, d} + \|H\|_{\theta, T, d} \big ).
\eeq
\end{theorem}
 \bpf The proof is  similar to the one  of Theorem \ref{t1}.
 Let $v \in {\cal C}^{2+ \theta}_{T, d} $ be a solution.
  One first proves the following maximum principle
 $$
 \sup_{t \in [0,T],\, x \in \R^n} | v  (t,x)|
  \; \le \;  T
   \sup_{t \in [0,T],\, x \in \R^n} | H(t,x)| +  \| g\|_0.
$$
 arguing as in \eqref{333} (using that
   $\aa = {\cal A}_1 +  {\cal A}_2$ as in the proof of Theorem \ref{t1}).

Concerning the localization procedure which gives the required a
priori estimates, we only   note that, for any $\eta \in
C_0^{\infty} (\R^n)$, according to the definition \eqref{cie},
 the function $v \eta$ solves
$$
 \left\{ \bal \partial_t \big( v \eta \big) (t,x)
  & =  \aa  ( \eta   v   ) (t,x)   -  v (t,x)\aa \eta
  (x)  - \lan Q(x) D \eta (x), D v(t,x) \ran
   +\, H(t,x) \eta(x),
 \; t \in (0,T],
 \\
      (\eta v)(0, x) & = \eta(x)g(x),\;\; x \in \R^n.
\eal \right.
$$
Finally the continuity method of Theorem \ref{t1} works also in
this case, replacing the space ${\cal C}^{2 + \theta}_d (\R^n)$
with
 ${\cal C}^{2 +  \theta}_{T,  d}$ and gives the
 assertion.
\epf \bre \label{fine} {\em One can weaken the assumption (ii) in
Hypothesis \ref{hyy}  about  $F$
  in order to prove
 elliptic and parabolic
  Schauder estimates for $\aa$.  To this purpose
  we can consider $F: \R^n \to \R^n$ such that
$F(x)= (F_1(x), \ldots, F_{\tilde p}(x), 0, \ldots, 0)$, $x \in
 \R^n$,  and moreover there exist $\theta \in (0,1)$ and
  $M>0$ such that,
  for any $x, y \in \R^n$, if $|y| \le 1$ then
 we have
 \beq \label{ne}
 |F(x) - F(x+ y) | \le M \nor  y \nor^{\theta}.
\eeq
 We  briefly explain how to prove  elliptic
  Schauder estimates for $\aa$ when $F$ satisfies the previous
 assumptions. First
    we deal with
 the maximum principle. Let $u \in {\cal C}
 ^{2+ \theta}_d (\R^n)$ be a solution. We consider $u_m =
  u * \rho_m$, where $(\rho_m)$ are mollifiers. Under the new
  assumptions on $F$ one can only show that $\aa u_m $ converges to $\aa u$
   uniformly on compact sets of $\R^n$ (compare with the proof
    of Theorem \ref{t1}). This fact allows
    to prove that if $x_0$ is a local maximum for $u$ then $\aa u
     (x_0) \le 0$ (see the proof of \cite[Proposition 3.1.10]{Lu2}).
      Adapting  the proof of \cite[Proposition 2.2]{Lu3}
  one obtains the maximum principle.
 Then, in order to get Schauder estimates, one writes
$$ \mbox{$
\lambda u (x) -
  \frac{1}{2} {\mbox{\rm Tr }}  (Q (x) D^2 u (x)) -
  {\langle}  Ax + (F * \rho) (x), D u (x) {\rangle}
    = f +
 {\langle}   F(x) - (F * \rho) (x) , D u (x){\rangle},
$}$$
 where $F * \rho$ is the convolution between $F$ and a
 function $\rho \in C^{\infty}_0 (\R^n)$,
   $\| \rho\|_{L^1(\R^n) } =1$, $0 \le \rho \le 1 $ and   $\rho (x)
=0$ if $|x| \ge 1$. Using that $D_i (F * \rho)(x) = $
  $ \int_{\R^n} (F(x-y) - F(x)) D_i \rho(y) dy$
 and similar formulae for higher
  partial derivatives,
   we see that $F * \rho$ satisfies  (iii) in Hypothesis \ref{hy}.
 Moreover by \eqref{ne} one checks that $F - (F*\rho)$ belongs to ${\cal
C}^{\theta}_d (\R^n)$. Straightforward computations allow to get
 Schauder estimates for $\aa$.
} \ere


\end{document}